\documentclass[fleqn,preprint,3p,times,sort&compress]{elsarticle}
 
\pdfoutput=1
\journal{}
\usepackage[hidelinks]{hyperref}
\hypersetup{colorlinks=true,linkcolor=blue,urlcolor=blue}
\usepackage{amsmath}
\usepackage{amsfonts}
\usepackage{amsbsy}
\usepackage{mathtools}
\usepackage{makecell}
\usepackage{booktabs}
\usepackage{multirow}
\usepackage[table,xcdraw,dvipsnames]{xcolor}
\usepackage{graphicx}
\usepackage{float}
\usepackage{overpic} 
\usepackage{rotating}
\usepackage{algorithm}
\usepackage{algpseudocode}
\usepackage{setspace}
\usepackage[labelfont=bf,labelsep=period]{caption}
\usepackage[labelfont=rm,font=footnotesize]{subcaption}
\numberwithin{equation}{section}

\makeatletter
\providecommand{\doi}[1]{%
  \begingroup
    \let\bibinfo\@secondoftwo
    \urlstyle{rm}%
    \href{http://dx.doi.org/#1}{%
      doi:\discretionary{}{}{}%
      \nolinkurl{#1}%
    }%
  \endgroup
}
\makeatother

\DeclareSymbolFont{rmlargesymbols}{OMX}{cmex}{m}{n}
\DeclareMathSymbol{\rmsum}{\mathop}{rmlargesymbols}{80}
\DeclareMathSymbol{\rmintop}{\mathop}{rmlargesymbols}{82}
\renewcommand{\sum}{\rmsum}
\renewcommand{\int}{\rmintop\nolimits}

\DeclareSymbolFont{greeksymbolsptm}{OML}{ptm}{m}{n}
\DeclareMathSymbol{\sigma}{\mathalpha}{greeksymbolsptm}{27}

\DeclareSymbolFont{greeksymbolsit}{OML}{mdbch}{i}{n}
\DeclareMathSymbol{\lambda}{\mathalpha}{greeksymbolsit}{21}
\DeclareMathSymbol{\nu}{\mathalpha}{greeksymbolsit}{23}

\DeclareSymbolFont{greeksymbols}{OML}{mdbch}{m}{n}
\DeclareMathSymbol{\Omega}{\mathalpha}{greeksymbols}{10}

\DeclareMathAlphabet{\altmathcal}{OMS}{cmsy}{m}{n}
\DeclareMathAlphabet{\altmathcalb}{OMS}{cmsy}{b}{n}
\newcommand{\Cr}{\altmathcal{C}}		
\newcommand{\Nr}{\altmathcal{N}}		
\newcommand{\Or}{\altmathcal{O}}		
\newcommand{\Sr}{\altmathcal{S}}		
\newcommand{\Lr}{\altmathcalb{L}}		
\newcommand{\Qr}{\altmathcalb{Q}}		
\newcommand{\Vr}{\altmathcalb{V}}		

\newcommand{\A}{\textbf{A}}				
\newcommand{\B}{\textbf{B}}				
\newcommand{\K}{\textbf{K}}				
\newcommand{\M}{\textbf{M}}				
\newcommand{\C}{\textbf{C}}				
\newcommand{\LL}{\textbf{L}}			
\renewcommand{\S}{\textbf{S}}			
\newcommand{\DD}{\textbf{D}}			
\renewcommand{\H}{\textbf{H}}			
\newcommand{\I}{\textbf{I}}				
\newcommand{\V}{\textbf{V}}				
\newcommand{\uu}{\textbf{u}}			
\newcommand{\uuh}{\uu_h}				
\newcommand{\tu}{\tilde{u}}     		
\newcommand{\tuu}{\tilde{\uu}}   		
\newcommand{\tuuh}{\tilde{\uu}_h}   	
\newcommand{\vv}{\textbf{v}}			
\newcommand{\e}{\textbf{e}}				
\newcommand{\rr}{\textbf{r}}			
\newcommand{\zz}{\textbf{z}}			
\newcommand{\hh}{\textbf{h}}			
\newcommand{\om}{\omega}				
\newcommand{\omh}{\om_h}				
\newcommand{\Om}{\Omega}				
\newcommand{\Gm}{\Gamma}				
\newcommand{\sh}{s}						
\newcommand{\lm}{\lambda}				
\newcommand{\lmh}{\lm_h}				
\newcommand{\eps}{\varepsilon}
\newcommand{\sg}{\sigma}

\newcommand{\E}{E}						
\newcommand{\EE}{\textbf{\E}}			
\newcommand{\EEh}{\EE_h}     			
\newcommand{\tE}{\tilde{\E}}     		
\newcommand{\tEE}{\tilde{\EE}_h}   		
\newcommand{\W}{W}						
\newcommand{\WW}{\textbf{\W}}			
\newcommand{\WWh}{\WW_h}     			
\newcommand{\w}{w}						
\newcommand{\ww}{\textbf{\w}}			
\newcommand{\wwh}{\ww_h}     			
\newcommand{\R}{\mathbb{R}}				
\newcommand{\CC}{\mathbb{C}}			
\newcommand{\ii}{{\rm i\.}}				
\newcommand{\bn}{\textbf{n}}			
\newcommand{\bF}{\textbf{F}}			
\renewcommand{\Re}{{\rm Re}\.}
\renewcommand{\Im}{{\rm Im}\.}

\newcommand{\curl}{{\rm curl}}
\renewcommand{\div}{{\rm div}}
\newcommand{\ibH}{\textit{\textbf{H}}}	
\newcommand{\LOM}{\big(L^2(\Om)\big)^2}

\newcommand{\n}{n} 						
\newcommand{\p}{p} 						
\newcommand{\Ba}[2]{B_{#1}^{\.#2}}		
\newcommand{\nx}{n_x} 					
\newcommand{\ny}{n_y} 					
\newcommand{\px}{p_x} 					
\newcommand{\py}{p_y} 					
\renewcommand{\ne}{n_e} 				
\newcommand{\N}{N} 				 		
\newcommand{\Nev}{N_{\rm eig}} 			
\newcommand{\Nit}{\Nr_{\rm Iter}}		
\newcommand{\Nfa}{\Nr_\fa}				
\newcommand{\Nfb}{\Nr_\fb}				
\newcommand{\Nmv}{\Nr_\mv}				
\newcommand{\Nvv}{\Nr_\Vv}				
\newcommand{\nnz}{N_{nz}}				
\newcommand{\fa}{{\rm Fact}}			
\newcommand{\fb}{{\rm FB}}				
\newcommand{\Vv}{{\rm VV}}				
\newcommand{\mv}{{\rm MV}}				
\newcommand{\IGA}{{\rm IGA}}		
\newcommand{\rIGA}{{\rm rIGA}}

\renewcommand{\.}{\!\:}									
\newcommand{\ceq}{\coloneqq}							
\newcommand{\pra}[1]{\left(#1\right)}					
\newcommand{\prc}[1]{\left\{#1\right\}} 				
\newcommand{\smaller}[1]{{\textstyle\footnotesize#1}} 	
\newcommand{\norm}[1]{\left\lVert#1\right\rVert} 		
\newcommand{\normB}[1]{\norm{#1}_{\.\B}}

\newcommand{\normL}[1]{\norm{#1}_{\.L^2}}

\newcommand{\az}[1]{{\color{blue}#1}}
\newcommand{\fig}[1]{\mbox{\az{Fig.} #1}}
\newcommand{\figs}[2]{\mbox{\az{Figs.} #1 and #2}}

\newcommand{\Eq}[1]{\mbox{Equation #1}}

\newcommand{\Eqs}[2]{\mbox{Equations #1 and #2}}

\newcommand{\tab}[1]{\mbox{\az{Table} #1}}

\newcommand{\Sec}[1]{\mbox{Section #1}}

\newcommand{\Rem}[1]{\mbox{\az{Remark} #1}}

\newcommand{\Alg}[1]{\mbox{\az{Algorithm} #1}}

\newdefinition{remark}{Remark}
\newdefinition{lem}{Lemma}

\definecolor{drot}{rgb}{0.7,0,0.1}

\begin{document}
\baselineskip14pt
\sloppy

\begin{frontmatter}
	
\title{Performance of Refined Isogeometric Analysis in Solving\\Quadratic Eigenvalue Problems}

\author[a1]{Ali Hashemian\corref{cor1}}
\cortext[cor1]{Corresponding author}
\ead{ahashemian@bcamath.org}

\author[a5]{Daniel Garcia}
\author[a2,a1,a3]{David Pardo}
\author[a6]{Victor M. Calo} 

\address[a1]{BCAM -- Basque Center for Applied Mathematics, Bilbao, Basque Country, Spain}
\address[a5]{IDAEA -- Institute of Environmental Assessment and Water Research, Barcelona, Catalu\~na, Spain}
\address[a2]{University of the Basque Country UPV/EHU, Leioa, Basque Country, Spain}
\address[a3]{Ikerbasque -- Basque Foundation for Sciences, Bilbao, Basque Country, Spain}
\address[a6]{Faculty of Science and Engineering, Curtin University, Perth, Australia}

\begin{abstract}
Certain applications that analyze damping effects require the solution of quadratic eigenvalue problems (QEPs). We use refined isogeometric analysis (rIGA) to solve quadratic eigenproblems. rIGA discretization, while conserving desirable properties of maximum-continuity isogeometric analysis (IGA), reduces the interconnection between degrees of freedom by adding low-continuity basis functions. This connectivity reduction in rIGA's algebraic system results in faster matrix LU factorizations when using multifrontal direct solvers.  We compare computational costs of rIGA versus those of IGA when employing Krylov eigensolvers to solve quadratic eigenproblems arising in 2D vector-valued multifield problems.  For large problem sizes, the eigencomputation cost is governed by the cost of LU~factorization, followed by costs of several matrix--vector and vector--vector multiplications, which correspond to Krylov projections.  We minimize the computational cost by introducing $C^0$ and $C^1$ separators at specific element interfaces for our rIGA generalizations of the curl-conforming N\'ed\'elec and divergence-conforming Raviart--Thomas finite elements.  Let $\p$ be the polynomial degree of basis functions; the LU factorization is up to ${\Or\big((\p-1)^2\big)}$ times faster when using rIGA compared to IGA in the asymptotic regime. Thus, rIGA theoretically improves the total eigencomputation cost by ${\Or\big((\p-1)^2\big)}$ for sufficiently large problem sizes. Yet, in practical cases of moderate-size eigenproblems, the improvement rate deteriorates as the number of computed eigenvalues increases because of multiple matrix--vector and vector--vector operations. Our numerical tests show that rIGA accelerates the solution of quadratic eigensystems by ${\Or(\p-1)}$ for moderately sized problems when we seek to compute a reasonable number of eigenvalues.
\end{abstract}

	
\end{frontmatter}


\section{Introduction} 
\label{sec.Introduction}

Refined isogeometric analysis (rIGA), introduced by~\citet{Garcia2017}, proved to be a successful extension of isogeometric analysis (IGA)~\cite{ Hughes2005} when using multifrontal direct solvers to approximate solutions to partial differential equations (c.f.,~\cite{ Garcia2019, PASZYNSKI2018, Siwik2019, Hashemian2021, Hashemian20212}).  rIGA preserves the desirable properties of maximum-continuity IGA discretizations while partitioning the computational domain into {\em macroelements} interconnected by low-continuity basis functions. This reduced connectivity exploits the recursive partitioning capability of multifrontal direct solvers, significantly reducing the solution cost.  In practice, when using rIGA discretizations for gradient-conforming $H^1$ spaces, the matrix factorization is asymptotically ${\Or(\p^2)}$ times faster in large domains compared to IGA (where $\p$ is the polynomial order of B-spline bases). In this context, rIGA improves the performance of discretizations by curl-conforming \ibH(curl) and divergence-conforming \ibH(div) spaces by up to ${\Or\big((\p-1)^2\big)}$. In all cases, rIGA also outperforms traditional finite element analysis (FEA) when considering a fixed number of elements and using the same polynomial order (see~\cite{ Garcia2017, Garcia2019}).

IGA is a widely used numerical technique to solve eigenvalue problems (see, e.g.,~\cite{ Cottrell2006, HUGHES2008, BUFFA2010, Hughes2014, Puzyrev2017, Hosseini2018, Mazza2019, Deng2019, Deng20192, GAO2020}). Recently, the application of rIGA in eigenanalysis has been also investigated in~\cite{ Puzyrev2018, Hashemian2021}.  However, most of the studied cases in the context of both maximum-continuity and refined IGA frameworks are limited to the (linear) {\em generalized} eigenvalue problems, which have the following form:
\begin{align}
(\A-\lm\B)\.\uu=\textbf{0}\,,
\end{align}
resulting in real eigenvalues when \B~is positive definite and \A~is Hermitian. Quadratic eigenvalue problems (QEPs), commonly expressed as
\begin{align}
(\K+\lm\C+\lm^2\M)\.\uu=\textbf{0}\,,
\end{align}
are the most applicable subset of polynomial eigenvalue problems. They are of great interest in multiple scientific and engineering case studies such as the vibration of flexible mechanisms~\cite{ SADLER1997}, fluid--structure interactions~\cite{ Olson1989}, constrained least-squares~\cite{ Sima2004}, electromagnetic wave propagation in conductive media~\cite{ Cooke2000}, and acoustic fluid in a cavity with absorbing boundaries~\cite{ Bermdez2000}.  In this work, we investigate the beneficial effect of using rIGA on eigencomputation cost of solving QEPs.
The arising quadratic eigensystems have complex eigenvalues in many practical occasions. Thus, as opposed to the generalized eigenproblems, some numerical algorithms (e.g., spectrum slicing, see \Sec{\ref{sec.EigsolutionAlgorithms}}) are inapplicable to quadratic cases. Additionally, the number of computed eigenvalues plays an important role in the overall eigencomputation cost of QEPs.  As a result, the implementation of rIGA discretizations shows different levels of improvement in linear (generalized) and quadratic eigenproblems (c.f.,~\cite{ Hashemian2021} for rIGA cost improvements in solving generalized eigenproblems).

When solving a quadratic eigenproblem, we may compute the eigensolution by {\em linearizing} the original eigensystem. Frequently used Krylov eigensolvers commonly perform an Arnoldi recurrence to project the equivalent linearized eigenproblem onto Krylov subspace\footnote{For some specific cases, quadratic eigensolvers may alternatively employ the nonsymmetric or pseudo-Lanczos algorithms (see, e.g,~\cite{ Campos2016, Campos2020}).}. Then, they find a finite number of eigenvalues using a Rayleigh--Ritz approximation (see, e.g.,~\cite{ Demmel2000, Tisseur2001, Campos2016, Campos2020}).  The projection process entails one (or more) LU factorization and multiple matrix--vector operations ---\.which are not naturally as expensive as matrix factorization but should be taken into account because the eigencomputation entails a significant number of these operations.
Eigencomputations are among the most expensive numerical approaches, especially for multidimensional systems and moderate to large problem sizes when we seek to compute many eigenvalues.  \citet{ Hashemian2021} identify the three most expensive numerical operations of solving generalized Hermitian eigenproblems, namely: matrix factorization, forward/backward eliminations (i.e., multiplications of LU factors by vectors), and matrix--vector products (i.e., multiplications of system matrices by vectors).  They show for eigenproblems in $H^1$ that, when using multifrontal direct solvers, rIGA improves matrix factorization asymptotically by ${\Or(\p^2)}$ for a sufficiently large number of degrees of freedom. The improvement factor of forward/backward eliminations is asymptotically up to $\Or(\p)$. At the same time, the higher number of nonzero entries in the rIGA matrices slightly degrades the matrix--vector product cost. These improvement/degradation ratios also hold for quadratic eigensystems. The difference, however, is in the number of times the eigensolver calls each operation.  In addition to the mentioned numerical operations, the vector--vector product plays an important role in eigencomputation cost of QEPs. Despite being considered as an inherently cheap numerical operation, it has an important contribution to the total cost when computing a large number of eigenvalues (see~\Sec{\ref{sec.CostEig}}).

Herein, we use rIGA discretizations to solve quadratic eigensystems of multidimensional problems with large numbers of degrees of freedom.  We consider two vector-valued {\em multifield} model problems arising in electromagnetic wave propagation in coductive media, and vibration of an acoustic fluid contained in a cavity with absorbing walls (see~\Sec{\ref{sec.ModelProblems}}).  We use higher-order B-spline generalizations of the curl-conforming N\'ed\'elec and divergence-conforming Raviart--Thomas finite elements (see, e.g.,~\cite{ BUFFA2010, Vazquez2010, EVANS2013, Buffa20102}) to discretize our 2D electromagnetic and vibroacoustic problems, respectively. These spaces allow us to obtain solution fields free of spurious eigenmodes (i.e., modes with no physical meaning). We briefly describe IGA and rIGA discretizations in~\Sec{\ref{sec.IGArIGA}}.
\Sec{\ref{sec.EigsolutionAlgorithms}} describes the eigensolution procedure applied to QEPs.
In \Sec{\ref{sec.CostEig}}, we theoretically estimate the eigencomputation costs based on the numerical operations employed by the quadratic eigensolver. We show that rIGA improves the efficiency of the most expensive numerical operations of the eigenanalysis. Our theoretical analysis shows that an improvement of ${\Or\big((\p-1)^2\big)}$ in eigencomputation cost is asymptotically possible when employing rIGA in multifield problems discretized by \ibH(curl) and \ibH(div) spaces. Indeed, for sufficiently large problems, the matrix factorization governs the solution cost. However, in practical moderate-size problems, the numerical tests of \Sec{\ref{sec.Results}} show that rIGA reduces the eigencomputation cost by a factor of approximately ${\Or(\p-1)}$ when computing a reasonable number of eigenvalues. Finally, we draw the main conclusions in \Sec{\ref{sec.Conclusions}}.


\section{Model problems} 
\label{sec.ModelProblems}

We consider the eigensolutions of two different quadratic eigenvalue problems.  The first eigenproblem arises in electromagnetic wave propagation; the second one arises from the analysis of an acoustic (i.e., inviscid, compressible, barotropic) fluid contained in a cavity with absorbing walls.


\subsection{Electromagnetic wave propagation in a conductive medium} 
\label{sub:Electromagnetic}

Let us define the \ibH(\curl)-conforming functional spaces on the computational domain ${\Om\subset\R^2}$ as
\begin{align}
\ibH(\curl;\Om)      &\ceq 
    \bigg\{\WW=\smaller{\begin{bmatrix}\W_x\\\W_y\end{bmatrix}}\in\LOM:\nabla\times\WW\in\LOM\bigg\} \,,\\
\ibH_{0}(\curl;\Om)  &\ceq 
    \bigg\{\WW\in\ibH(\curl;\Om):\WW\times\bn=\textbf{0}{\rm~on~}\partial\Om\bigg\} \,,
\label{eq.HcurlSpace}
\end{align}
where $\bn$ is the outward unit normal vector on the boundary ${\partial\Om}$.
We consider the eigenvalue problem of the electromagnetic wave propagation equation in a {\em conductive} medium:
\begin{align}
\left\{\begin{matrix*}[l]
  {\rm Find~} \om\in\CC {\rm~and~} \EE=\smaller{\begin{bmatrix}\E_x\\\E_y\end{bmatrix}} \in\ibH_0(\curl;\Om),
  {\rm~with~} \EE:\Om\rightarrow\CC^2, {\rm~such~that}\\[7pt]
  \qquad\begin{aligned}
    \nabla\times(\,\tilde{\mu}^{-1}\nabla\times\EE) + \tilde{\pmb{\sg}}\EE &=\textbf{0}\,, 
      \quad&&{\rm in~}\Om \,,\\[3pt]
    \EE\times\bn &= \textbf{0}\,, 
      \quad&&{\rm on~}\partial{\Om} \,,\end{aligned}
\end{matrix*}\right.
\label{eq.HcurlStrong1}
\end{align}
where $\om$ is the eigenfrequency, $\EE$ is the electromagnetic field, 
${\tilde{\pmb{\sg}}\ceq\pmb{\sg}+\ii\om\eps\.\I}$ and ${\tilde{\mu}\ceq\ii\om\mu}$
with $\pmb{\sg}$ being the conductivity matrix, $\ii$ being the imaginary unit, and $\eps$ and $\mu$ being the electric permittivity and magnetic permeability, respectively. 
We rewrite the strong form of~\eqref{eq.HcurlStrong1} as
\begin{align}
\nabla\times(-\mu^{-1}\nabla\times\EE)- \ii\om\.\pmb{\sg}\EE+\eps\om^2\EE=\textbf{0}\,,
\label{eq.HcurlStrong2}
\end{align}
and consider ${\WW\in\ibH_{0}(\curl;\Om)}$ as an arbitrary test function to construct the sesquilinear forms
\begin{align}
a\.(\WW,\EE)&\ceq\int_\Om (\overline{\nabla\times\WW})\cdot(\nabla\times-\mu^{-1}\EE)\,d\Om\,,\\
b\.(\WW,\EE)&\ceq\int_\Om \overline{\WW}\cdot\eps\.\EE\,d\Om\,,\\
c\.(\WW,\EE)&\ceq\int_\Om \overline{\WW}\cdot(-\.\ii\pmb{\sg}\EE)\,d\Om\,,
\label{eq.HcurlQuadMatrices}
\end{align}
where {\em overbar} means the conjugate form of the complex vector space and $\cdot$ denotes the inner product. We build the weak form of~\eqref{eq.HcurlStrong1} as follows:
\begin{align}
a\.(\WW,\EE)+\om\.c\.(\WW,\EE)+\om^2b\.(\WW,\EE)=0\,.
\label{eq.HcurlQuadWeak}
\end{align}

Let us define
\begin{align}
\Vr_{0,h}^{\.\curl}(\Om)\ceq\prc{\WWh\in\ibH(\curl;\Om)\cap\ibH^1(\Om):\WWh\times\bn=\textbf{0}{\rm~on~}\partial{\Om}},
\label{eq.HcurlDiscreteSpace}
\end{align}
as our curl-conforming approximate space, where index $h$ refers to approximate values obtained by discretization.
We use a Galerkin discretization of the continuous eigenproblem and express the approximate solution field as follows:
\begin{align}
\EE_h=\sum_{i=0}^{\N-1}\pmb{\phi}_i^{\.\curl}\tE_{h,i}\,,
\label{eq.HcurlApproxSolution}
\end{align}
where $\tE_{h,i}$ represents the ${i\.}$-th nodal component of the discrete solution $\tEE$ given by 
\begin{align}
\tEE=\big[\tE_{x,0}\.,...\.,\tE_{x,\N_x-1},\tE_{y,0}\.,...\.,\tE_{y,\N_y-1}\big]\,,
\label{eq.HcurlControlVariables}
\end{align}
and ${\pmb{\phi}_i^{\.\curl}}$ is its corresponding curl-conforming vector-valued basis function. Herein, ${\N\ceq\N_x+\N_y}$ is the total number of degrees of freedom, being $\N_x$ and $\N_y$ the numbers of degrees of freedom of vector fields in $x$ and $y$ directions, respectively.
Thus, we write the approximate form of~\eqref{eq.HcurlQuadWeak} as
\begin{align}
a\.(\WWh\.,\EEh)+\om_h\.c\.(\WWh\.,\EEh)+\om_h^2\.b\.(\WWh\.,\EEh)=0\,,
\label{eq.HcurlQuadWeakDiscrete}
\end{align}
that results in the following quadratic eigenvalue problem in the discrete form:
\begin{align}
\big(\K+\omh\C+\omh^2\.\M\big)\.\tEE=\textbf{0} \,.
\label{eq.HcurlQuadEig}
\end{align}

\begin{remark}
  In the context of eigenvalue analysis, we refer to the discrete solution vector $\tEE$ as the {\em eigenvector} and the continuous solution filed (either exact $\EE$ or approximate $\EEh$) as the {\em eigenfunction}. For an $N$-degree-of-freedom system, the eigenproblem~\eqref{eq.HcurlQuadEig} has $2N$ approximate {\em eigenpairs} referring to either ${(\omh,\EEh)}$ or ${(\omh,\tEE)}$\..
\end{remark}

The {\em stiffness} and {\em mass} matrices in eigenproblem~\eqref{eq.HcurlQuadEig}, $\K$ and $\M$, respectively, are symmetric and $\M$ is positive definite. The {\em damping} matrix $\C$ is skew-Hermitian (i.e., ${\C=-\C^T}$, where superscript $T$ denotes the conjugate transpose). 
Thus, the current model problem results in a {\em gyroscopic} eigensystem (c.f.,~\cite{Qian2007}) with a Hamiltonian spectrum structure whose eigenvalues are symmetric with respect to the imaginary axis.
We refer to ${\Re(\omh)}$ and ${\Im(\omh)\geq0}$ as the angular frequency and exponential decay of the electromagnetic wave with time, respectively~\cite{Cooke2000}.

\begin{remark}
  In a non-conductive medium with ${\pmb{\sg}=\textbf{0}}$\., the eigensystem~\eqref{eq.HcurlQuadEig} is reduced to the generalized eigenproblem ${(\A-\lmh\B)\.\tEE=\textbf{0}}$ and commonly referred to as the Maxwell's eigenvalue problem (i.e., the curl-curl operator eigenproblem~\cite{ BUFFA2010, Zhang2018}).  We set ${\lmh\ceq\omh^2}$ and define ${\A\ceq-\K}$ and ${\B\ceq\M}$ to follow the notations of generalized eigenproblems described in \Sec{\ref{sec.EigsolutionAlgorithms}}. In this case, since $\B$ is positive definite and $\A$ is positive semidefinite, all $\lmh$ values are real non-negative and the electromagnetic wave does not decay.
\label{rem.LinearHcurl}
\end{remark}


\subsection{Acoustic fluid contained in a cavity with absorbing walls} 
\label{sub:AcousticFluid}

We consider an acoustic fluid contained in a rigid cavity with some (or all) of its walls covered by a viscoelastic material absorbing the acoustic energy of the fluid~\cite{Bermdez2000}.
Let us consider the domain ${\Om\subset\R^2}$ and its boundary ${\partial\Om\ceq\Gm_A\cup\Gm_R}$\., being $\Gm_A$ and $\Gm_R$ the union of {\em absorbing} and {\em rigid} boundaries, respectively.
We define the \ibH(\div)-conforming functional spaces as follows:
\begin{align}
\ibH(\div;\Om)      &\ceq 
  \bigg\{\ww=\smaller{\begin{bmatrix}\w_x\\\w_y\end{bmatrix}}\in\LOM:\nabla\cdot\ww\in L^2(\Om)\bigg\} \,,\\
\ibH_{0}(\div;\Om)  &\ceq 
  \bigg\{\ww\in\ibH(\div;\Om):\ww\cdot\bn\in L^2(\partial\Om){\rm~and~}\ww\cdot\bn=0~{\rm on }~\Gm_R\bigg\}\,.
\label{eq.HdivSpace}
\end{align}
We find the damped vibration modes of the fluid as complex solutions to the following problem:
\begin{align}
\left\{\begin{matrix*}[l]
  {\rm Find~} \lm\in\CC, \uu=\smaller{\begin{bmatrix}u_x\\u_y\end{bmatrix}}\in \ibH_0(\div;\Om)
  {\rm ~and~} P\in H_0^1(\Om),{\rm~with~} 
  \uu:\Om\rightarrow\CC^2 {\rm ~and~} P:\Om\rightarrow\CC,
  {\rm~such~that}\\[7pt]
  \qquad\begin{aligned}
    \rho\lm^2\.\uu+\nabla P			    &=0\,, \quad&&{\rm in~}\Om \,,\\[2pt]
   -\rho c^2\.\nabla\cdot\uu		    &=P\,, \quad&&{\rm in~}\Om \,,\\[2pt]
    (\alpha+\lm\.\beta)\.\uu\cdot\bn&=P\,, \quad&&{\rm on~}\Gm_A \,,\\[2pt]
    \uu\cdot\bn 					          &=0\,, \quad&&{\rm on~}\Gm_R \,,\end{aligned}
\end{matrix*}\right.
\label{eq.HdivStrong}
\end{align}
where 
$P$ is the fluid pressure,
$\uu$ is the displacement field,
$\rho$ is the fluid density, and
$c$ is the acoustic speed.
Herein, ${\Im(\lm)}$ and ${\Re(\lm)}$ are the vibration's angular frequency and decay rate, respectively, while coefficients $\alpha$ and $\beta$ are related to the impedance of the viscoelastic material.

By considering an arbitrary test function ${\ww\in\ibH_0(\div;\Om)}$ and defining the sesquilinear forms
\begin{align}
\breve{a}\.(\ww,\uu)&\ceq\int_\Om (\nabla\cdot\overline{\ww})\.(\,\rho c^2\nabla\cdot\uu)\,d\Om+
\int_{\Gm_A} (\overline{\ww}\cdot\bn)\.(\alpha\.\uu\cdot\bn)\,d\.\Gm_A\,,
\label{eq.HdivQuadMatrces.K}\\
\breve{b}\.(\ww,\uu)&\ceq\int_\Om \overline{\ww}\cdot\rho\.\uu\,d\Om\,,
\label{eq.HdivQuadMatrces.C}\\
\breve{c}\.(\ww,\uu)&\ceq\int_{\Gm_A} (\overline{\ww}\cdot\bn)\.(\,\beta\.\uu\cdot\bn)\,d\.\Gm_A\,,
\label{eq.HdivQuadMatrces.M}
\end{align}
we obtain the weak form of~\eqref{eq.HdivStrong} involving only displacement variables as follows (we use the ${\,\breve{}\,}$ symbol to distinguish between sesquilinear forms of our two model problems):
\begin{align}
\breve{a}\.(\ww,\uu)+\lm\.\breve{c}\.(\ww,\uu)+\lm^2\.\breve{b}\.(\ww,\uu)=0\,.
\label{eq.HdivWeak}
\end{align}
We write another variational formulation in terms of the pressure that leads to a nonlinear eigenproblem~\cite{ Bermdez2000} that is not studied here.  Let us consider
\begin{align}
\Vr_{0,h}^{\.\div}(\Om)\ceq\prc{\wwh\in\ibH(\div;\Om)\cap\ibH^1(\Om):\wwh\cdot\bn\in L^2(\partial\Om){\rm~and~}\wwh\cdot\bn=0{\rm~on~}\Gm_R},
\label{eq.HdivDiscreteSpace}
\end{align}
as our divergence-conforming approximate space. We discretize the continuous displacement field as
\begin{align}
\uuh=\sum_{i=0}^{\N-1}\pmb{\phi}_i^{\.\div}\.\tu_{h,i}\,,
\label{eq.HdivApproxSolution}
\end{align}
where ${\pmb{\phi}_i^{\.\div}}$ is the divergence-conforming vector-valued basis and $\tu_{h,i}$ is its corresponding component of the discrete solution field ${\tuuh=\big[\tu_{x,0}\.,...\.,\tu_{x,\N_x-1},\tu_{y,0}\.,...\.,\tu_{y,\N_y-1}\big]}$. Thus, we write the approximate form of~\eqref{eq.HdivWeak} as
\begin{align}
\breve{a}\.(\wwh\.,\uuh)+\lmh\.\breve{c}\.(\wwh\.,\uuh)+\lmh^2\.\breve{b}\.(\wwh\.,\uuh)=0\,,
\label{eq.HdivWeakDiscrete}
\end{align}
resulting in 
\begin{align}
\big(\breve{\K}+\lmh\breve{\C}+\lmh^2\.\breve{\M}\big)\.\tuuh=\textbf{0} \,,
\label{eq.HdivQuadEig}
\end{align}
which is a Hermitian quadratic eigenproblem (all system matrices are real symmetric). The mass matrix~$\breve{\M}$ is positive definite, while the stiffness $\breve{\K}$ and damping $\breve{\C}$ are positive semidefinite matrices resulting in a stable vibration response, i.e., ${\Re(\lmh)\leq0}$.
We compute the eigenvalues as a combination of purely real and conjugate complex pairs.  The former corresponds to an overdamped non-oscillatory response, while the latter amounts to an underdamped vibration.


\section{Refined isogeometric analysis} 
\label{sec.IGArIGA}

We now consider multifield vectorial solutions and discretize the continuous eigenproblems~\eqref{eq.HcurlStrong1} and~\eqref{eq.HdivStrong} using B-spline generalizations of curl-conforming N\'ed\'elec and divergence-conforming Raviart--Thomas spaces, respectively (see, e.g.,~\cite{Garcia2019}).
We herein review some basic concepts of maximum-continuity and refined IGA discretizations.


\subsection{Bivariate B-spline space} 
\label{sub:Bsplines}

Given the parametric domain ${\hat{\Om}:\xi,\eta\in(0,1)^2\subset\R^{2}}$,
we introduce the bivariate B-spline space ${\Sr^{\px,\py}_{k_x,k_y}}$ as
\begin{align}
\Sr^{\px,\py}_{k_x,k_y}\ceq{\rm span}\prc{\Ba{i,j}{\px,\py}}_{i=0,j=0}^{\nx-1,\ny-1},
\end{align}
where $\n$, $\p$, and $k$ with their indices are the number of degrees of freedom, polynomial degree, and continuity of basis functions in each direction, respectively. The bivariate basis functions ${\Ba{i,j}{\px,\py}}$ are given in the tensor product sense by:
\begin{align}
\Ba{i,j}{\px,\py}\ceq\Ba{i}{\px}(\xi)\.\Ba{j}{\py}(\eta) \,, \quad i=0,1,...,\nx-1 \,, \quad j=0,1,...,\ny-1 \,,
\end{align}
where the univariate bases span their respective knot sequences
\begin{align}
\Xi & =[\underbrace{0,0,...,0}_{\px+1},\xi_{\px+1},\xi_{\px+2}...,\xi_{\nx-1},\underbrace{1,1,...,1}_{\px+1}]\,, \label{eq.KnotsX}\\
\rm{H} & =[\underbrace{0,0,...,0}_{\py+1},\eta_{\py+1},\eta_{\py+2}...,\eta_{\ny-1},\underbrace{1,1,...,1}_{\py+1}]\,. \label{eq.KnotsY}
\end{align}
We use the Cox--De Boor recursion formula~\cite{ TheNURBSBook} to obtain the univariate basis functions, which in terms of $\xi$ reads:
\begin{align}
\Ba{i}{0}(\xi)  &= 
\begin{cases}\begin{aligned}
&1\,,&&\xi_i\leq\xi<\xi_{i+1}\,, \\ \label{eq.BasisFunBi0} 
&0\,,&&{\rm otherwise}\,,
\end{aligned}\end{cases} \\
\Ba{i}{\px}(\xi) &= \dfrac{\xi-\xi_i}{\xi_{i+\px}-\xi_i}\Ba{i}{\px-1}(\xi)+
\dfrac{\xi_{i+\px+1}-\xi}{\xi_{i+\px+1}-\xi_{i+1}}\Ba{i+1}{\px-1}(\xi)\,. \label{eq.BasisFunBip}
\end{align}
When evaluating basis functions, we find the corresponding nonzero knot span in~\eqref{eq.BasisFunBi0} and efficiently evaluate~\eqref{eq.BasisFunBip} avoiding any division by zero. More details are given in~\cite[Algorithms A2.1 and A2.2]{TheNURBSBook}.


\subsection{IGA discretization} 
\label{sub:IGA}

In the IGA context, we consider single multiplicities for all {\em interior} knots in~\eqref{eq.KnotsX}~and~\eqref{eq.KnotsY}, providing maximum continuity of basis functions. Let us consider polynomial degree $\p$ for discretization and assume the 2D domain has an ${\ne\times\ne}$ grid, being $\ne$ the number of elements in each direction. We construct the curl- and divergence-conforming discrete spaces in the parametric domain (c.f.,~\cite{ BUFFA2010, Garcia2019}) as
\begin{align}
\hat{\Vr}_{h,\IGA}^{\,\curl}(\hat{\Om}) &\ceq \Sr^{\p-1,\p}_{k-1,k} \times \Sr^{\p,\p-1}_{k,k-1} \,, \\[3pt]
\hat{\Vr}_{h,\IGA}^{\,\div}(\hat{\Om}) &\ceq \Sr^{\p,\p-1}_{k,k-1} \times \Sr^{\p-1,\p}_{k-1,k} \,,
\label{eq.DiscreteSpacesIGA}
\end{align}
where ${k=\p-1}$. The following vector-valued bivariate B-spline bases characterize these spaces:
\begin{align}
\hat{\pmb{\phi}}^{\.\curl}_{\.\IGA} &\ceq 
\begin{bmatrix} \Ba{0,0}{\p-1,\p} &...& \Ba{\ne+p-2,\.\ne+p-1}{\p-1,\p} & 0 &...&0\\[3pt]
0 &...& 0 & \Ba{0,0}{\p,\p-1} &...& \Ba{\ne+p-1,\.\ne+p-2}{\p,\p-1}\end{bmatrix} \., \label{eq.HcurlBasis} \\[5pt]
\hat{\pmb{\phi}}^{\.\div}_{\.\IGA}  &\ceq 
\begin{bmatrix} \Ba{0,0}{\p,\p-1} &...& \Ba{\ne+p-1,\.\ne+p-2}{\p,\p-1} & 0 &...&0\\[3pt]
0 &...& 0 & \Ba{0,0}{\p-1,\p} &...& \Ba{\ne+p-2,\.\ne+p-1}{\p-1,\p}\end{bmatrix} \..
\label{eq.HdivBasis}
\end{align}

\begin{figure}[!h]
	\begin{subfigure}{0.49\textwidth} \hspace{-8pt}
		\includegraphics{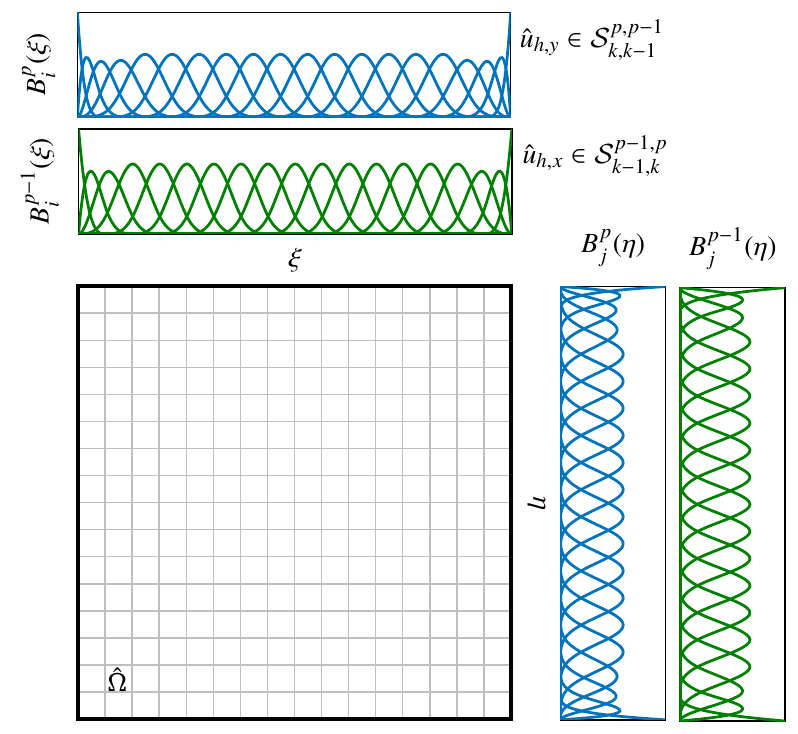}
		\caption{\ibH(\curl) IGA discrete space}
	\end{subfigure}
	\begin{subfigure}{0.49\textwidth} \hspace{6pt}
		\includegraphics{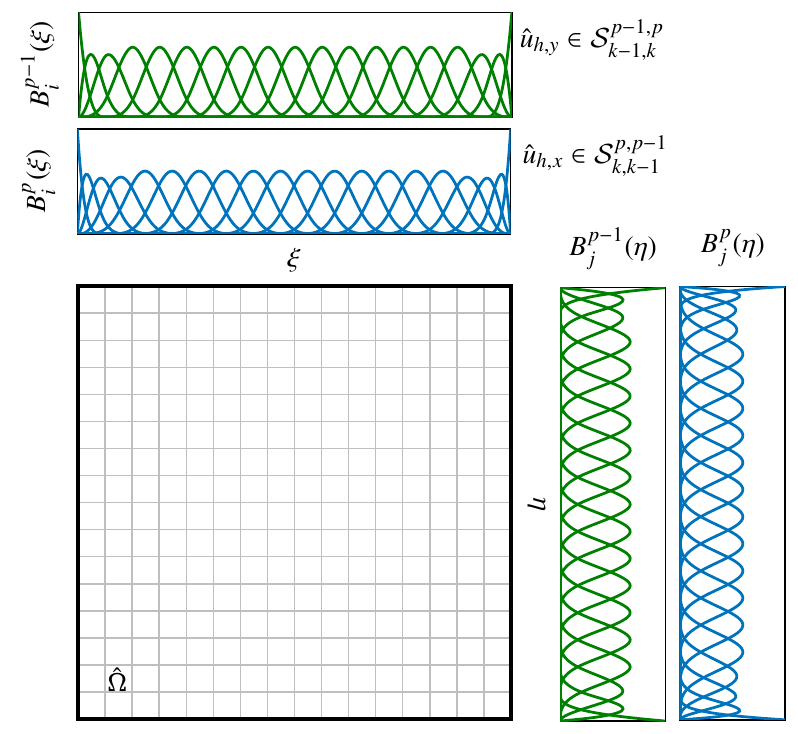}
		\caption{\ibH(\div) IGA discrete space}
	\end{subfigure}
	\caption{Examples of \ibH(\curl) and \ibH(\div) spaces discretized by maximum-continuity IGA with uniform ${16\times16}$ elements, polynomial degree ${\p=4}$, and continuity ${k=3}$. We show the univariate bases of the solution fields $\hat{u}_{h,x}$ and $\hat{u}_{h,y}$ in the parametric domain, where the blue and green color codes refer to basis functions of degree $\p$ and ${\p-1}$, respectively.}
	\label{fig.IGAHcurlHdiv}
\end{figure}

\fig{\ref{fig.IGAHcurlHdiv}} illustrates examples of \ibH(\curl) and \ibH(\div) IGA discrete spaces as well as univariate basis functions of the respective vector fields in the parametric domain. For 2D problems, \ibH(\curl) and \ibH(\div) are the same spaces rotated by 90 degrees in the parametric domain.  Herein, for the sake of simplicity, we use $\uu$ in reference to the solution fields of both model problems described in \Sec{\ref{sec.ModelProblems}}.

Let us define ${\bF:\hat{\Om}\rightarrow\Om}$ as the geometric mapping from the parametric space onto the physical domain and $D\bF$ as its Jacobian.
We introduce our discrete spaces in the physical domain as
\begin{align}
\Vr_{h,\IGA}^{\,\curl}(\Om) &\ceq \prc{\uuh\in\ibH(\curl;\Om)\cap\ibH^1(\Om):\iota^{\.\curl}(\uuh)=\hat{\uu}_h\in\hat{\Vr}_{h,\IGA}^{\,\curl}(\hat{\Om})} \,, \\
\Vr_{h,\IGA}^{\,\div}(\Om) &\ceq \prc{\uuh\in\ibH(\div;\Om)\cap\ibH^1(\Om):\iota^{\.\div}(\uuh)=\hat{\uu}_h\in\hat{\Vr}_{h,\IGA}^{\,\div}(\hat{\Om})} \,, 
\end{align}
considering the following curl- and divergence-preserving mappings~\cite{Garcia2019,BUFFA2010}:
\begin{align}
\iota^{\.\curl}(\uuh)	&\ceq (D\bF)^{T}(\uuh \circ \bF) \,, \label{eq.HcurlMapping}\\
\iota^{\.\div}(\uuh)	&\ceq \det\.(D\bF)\.(D\bF)^{-1}(\uuh \circ \bF)\,. \label{eq.HdivMapping}
\end{align}


\subsection{rIGA discretization} 
\label{sub:rIGA}

Refined isogeometric analysis (rIGA) is a discretization technique that optimizes the performance of direct solvers.  In particular, rIGA preserves the optimal convergence order of direct solvers for a fixed number of elements in the domain.  \citet{Garcia2017} first presented this strategy for $H^1$ spaces and then extended it to \ibH(\curl), \ibH(\div), and $L^2$ spaces~\cite{ Garcia2019}.  Starting from the maximum-continuity IGA discretization, rIGA reduces the continuity of certain basis functions by increasing the multiplicity of the respective existing knots.  Hence, we subdivide the computational domain into high-continuity macroelements interconnected by low-continuity hyperplanes.  These hyperplanes coincide with the locations of the {\em separators} at different partitioning levels of the multifrontal direct solvers.  Thus, when solving systems of equations arising in $H^1$-discretized problems, rIGA performs matrix factorization asymptotically ${\Or(\p^2)}$ times faster in large domains ---\.and $\Or(\p)$ faster in small domains\.--- compared to IGA. The improvement factors are ${\Or\big((\p-1)^2\big)}$ and ${\Or(\p-1)}$ for large and small problems, respectively, when considering \ibH(curl) and \ibH(div) spaces (mainly because we reduce to degree ${\p-1}$ at certain spatial directions).  In comparison to traditional FEA with the same number of elements, rIGA provides even larger improvements.  Additionally, rIGA reduces the forward/backward elimination cost asymptotically by up to $\Or(\p)$ in $H^1$ and ${\Or(\p-1)}$ in \ibH(curl) and \ibH(div) for sufficiently large domains since the LU (or Cholesky) factors have fewer nonzero terms for rIGA.  Consider a mesh with a fixed element number in each direction being a power of two; these improvements correspond to macroelements with~16 elements in each direction for moderately sized problems (see~\cite{ Garcia2017, Garcia2019, Hashemian2021}).

\begin{figure}[!h]
	\begin{subfigure}{0.49\textwidth} \hspace{-8pt}
		\includegraphics{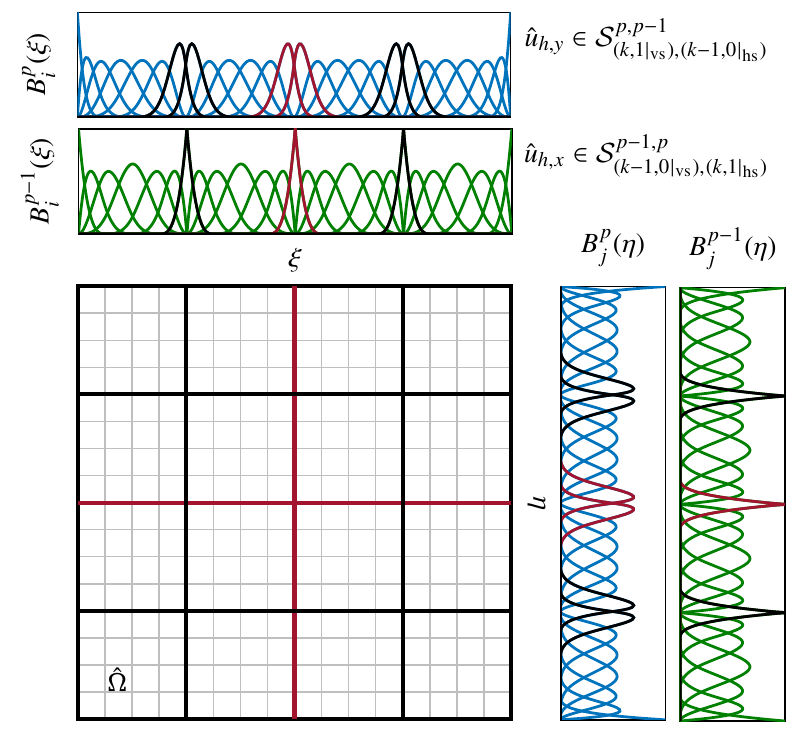}
		\caption{\ibH(\curl) rIGA discrete space}
	\end{subfigure}
	\begin{subfigure}{0.49\textwidth} \hspace{6pt}
		\includegraphics{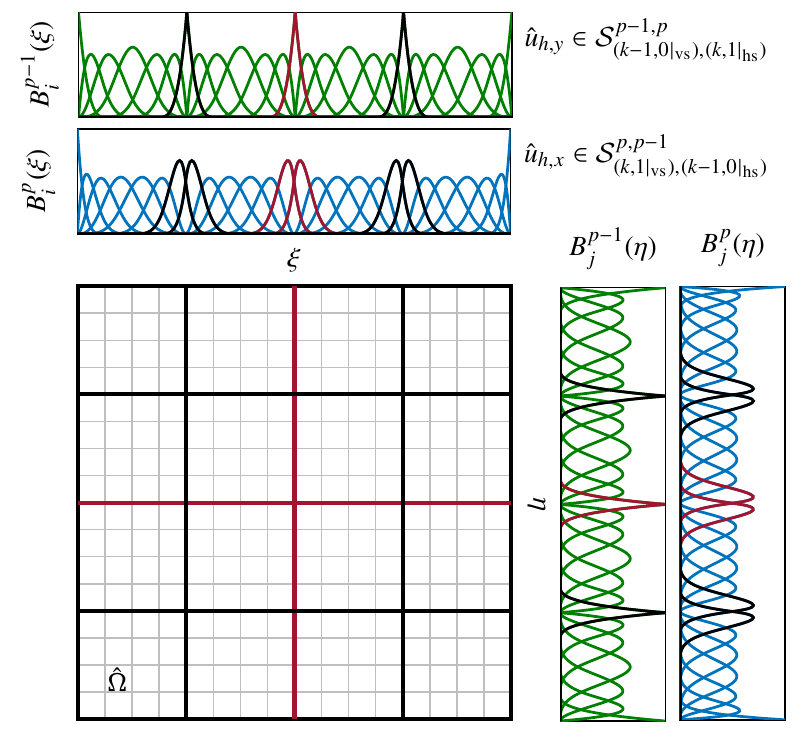}
		\caption{\ibH(\div) rIGA discrete space}
	\end{subfigure}
  \caption{\ibH(\curl) and \ibH(\div) rIGA spaces of the ${16\times16}$ domain of \fig{\ref{fig.IGAHcurlHdiv}} with ${\p=4}$ and ${k=3}$, after two levels of partitioning that generate four blocks in each direction. rIGA reduces the continuity of basis functions by ${k-1}$ degrees across the macroelement separators (we show the low-continuity bases corresponding to the first and second partitioning levels in dark red and black, respectively). Thin lines in the mesh skeleton denote the high-continuity element interfaces, while thick lines illustrate the macroelement boundaries.\vspace{4pt}}
	\label{fig.rIGAHcurlHdiv}
\end{figure}

We discretize our multifield problems using \ibH(\curl) and \ibH(\div) spaces; thus, we preserve the commutativity of the {\em de~Rham} diagram~\cite{ Demkowicz2000} by reducing the continuity in ${k-1}$ degrees.  We achieve the commutativity using both $C^{\.0}$ and $C^{1}$ hyperplanes and reduce the continuity across interfaces between subdomains (i.e., macroelements). Thus, we represent the rIGA spaces in the parametric domain as follows:
\begin{align}
\hat{\Vr}_{h,\rIGA}^{\,\curl}(\hat{\Om}) &\coloneqq \Sr^{p-1,p}_{(k-1,0\.|_{\.\rm vs}),\.(k,1|_{\.\rm hs})} \times \Sr^{p,p-1}_{(k,1|_{\.\rm vs}),\.(k-1,0\.|_{\.\rm hs})} \,, \\[3pt]
\hat{\Vr}_{h,\rIGA}^{\,\div}(\hat{\Om}) &\coloneqq \Sr^{p,p-1}_{(k,1|_{\.\rm vs}),\.(k-1,0\.|_{\.\rm hs})} \times \Sr^{p-1,p}_{(k-1,0\.|_{\.\rm vs}),\.(k,1|_{\.\rm hs})} \,,
\label{eq.DiscreteSpacesrIGA}
\end{align}
where ``vs" and ``hs" denote the vertical and horizontal separators, respectively. Accordingly, one updates~{\eqref{eq.HcurlBasis}--\eqref{eq.HdivMapping}} in a similar manner to write rIGA basis functions, mappings, and discrete spaces in the physical domain.
\fig{\ref{fig.rIGAHcurlHdiv}} depicts an example of an rIGA discretization of the discrete spaces of \fig{\ref{fig.IGAHcurlHdiv}} obtained by two levels of symmetric partitioning, which results in ${4\times4}$ macroelements.

\begin{figure}[!h]
	\begin{subfigure}{0.33\textwidth} 
		\begin{overpic}[width=1\linewidth,trim={0cm 0cm 7.25cm 0cm},clip] {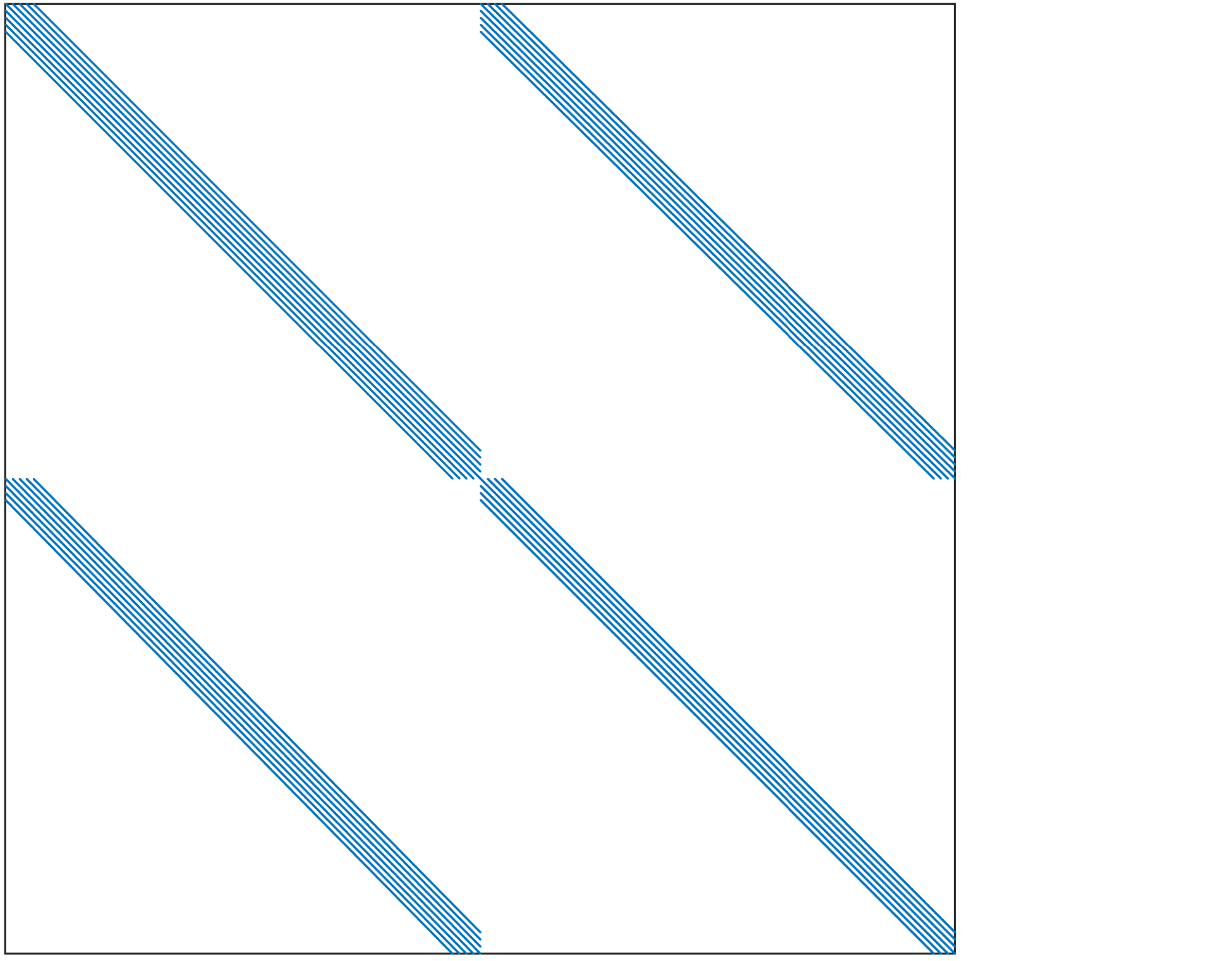}\end{overpic}
		\caption{IGA matrix: 9,112\,$\times$\,9,112}
		\label{fig.MatPatterns.IGAmatrix}
	\end{subfigure}
	\begin{subfigure}{0.33\textwidth} 
		\begin{overpic}[width=1\linewidth,trim={0cm 0cm 7.25cm 0cm},clip] {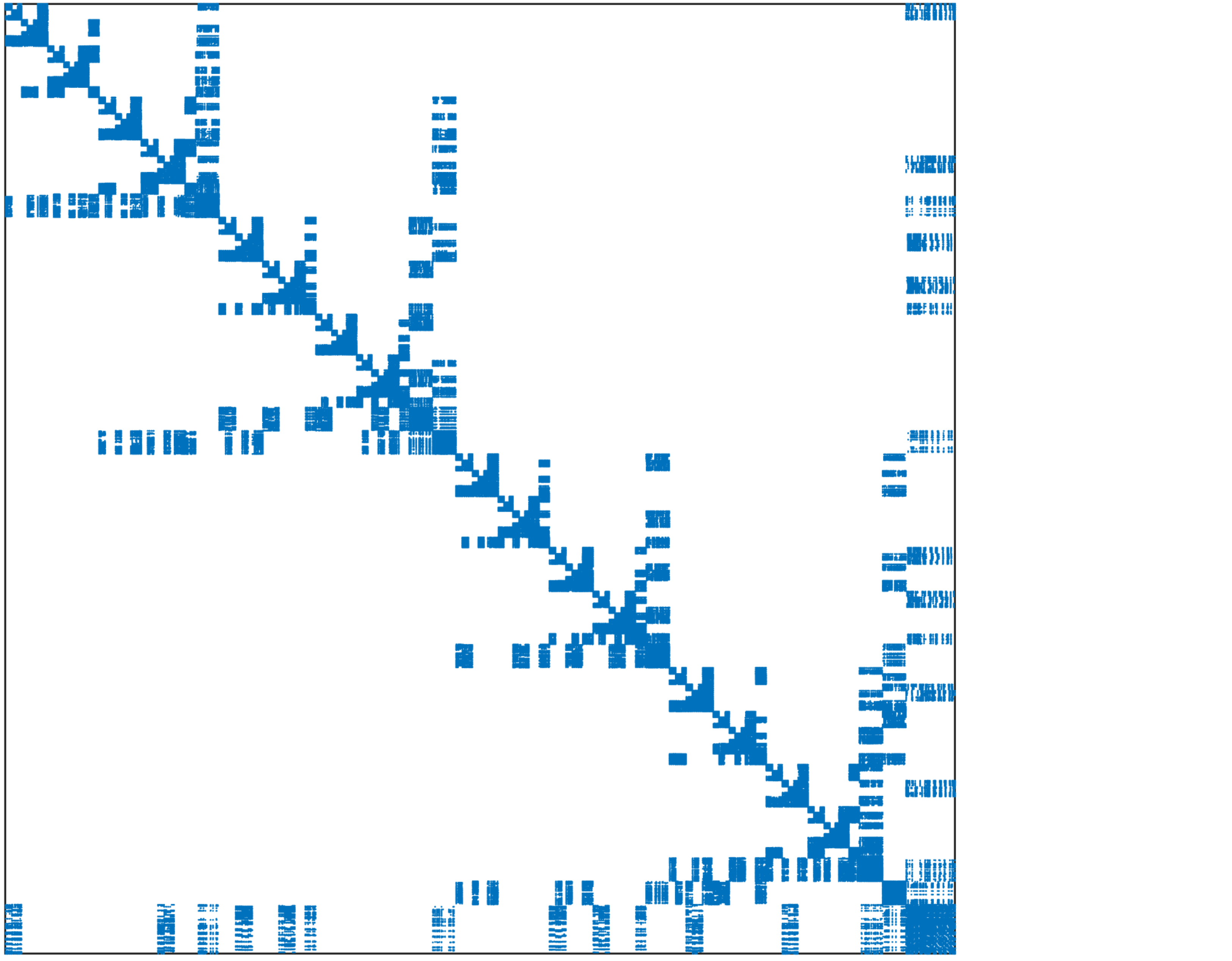}\end{overpic}
		\caption{IGA reordered matrix: $\nnz=\,$1,090,240}
		\label{fig.MatPatterns.IGAreorder}
	\end{subfigure}
	\begin{subfigure}{0.33\textwidth} 
		\begin{overpic}[width=1\linewidth,trim={0cm 0cm 7.25cm 0cm},clip] {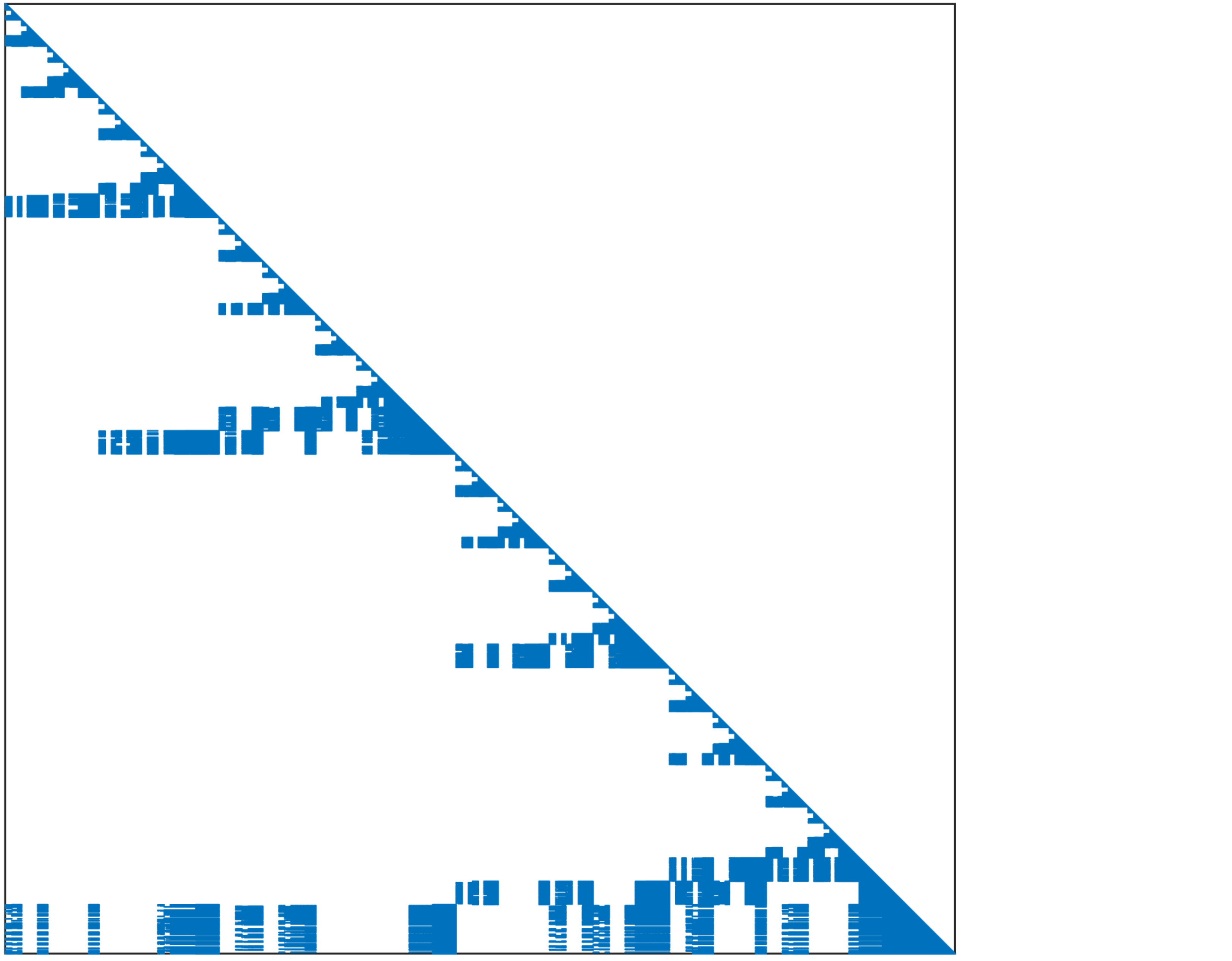}\end{overpic}
		\caption{IGA lower factor: $\nnz=\,$3,142,952}
		\label{fig.MatPatterns.IGAfactor}
	\end{subfigure}\\[5pt]
	\begin{subfigure}{0.33\textwidth} 
		\begin{overpic}[width=1\linewidth,trim={0cm 0cm 7.25cm 0cm},clip] {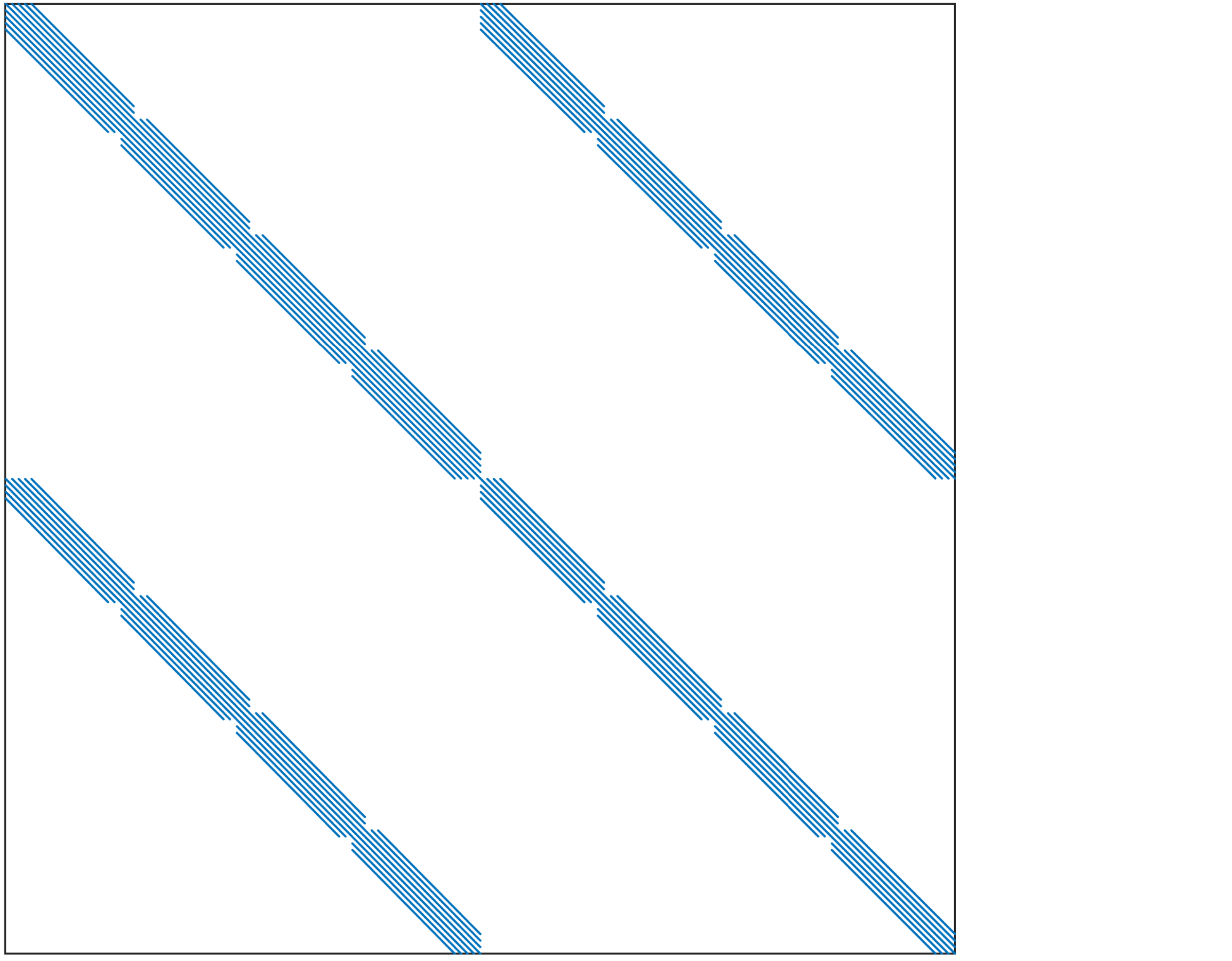}\end{overpic}
		\caption{rIGA matrix: 10,804\,$\times$\,10,804}
		\label{fig.MatPatterns.rIGAmatrix}
	\end{subfigure}
	\begin{subfigure}{0.33\textwidth} 
		\begin{overpic}[width=1\linewidth,trim={0cm 0cm 7.25cm 0cm},clip] {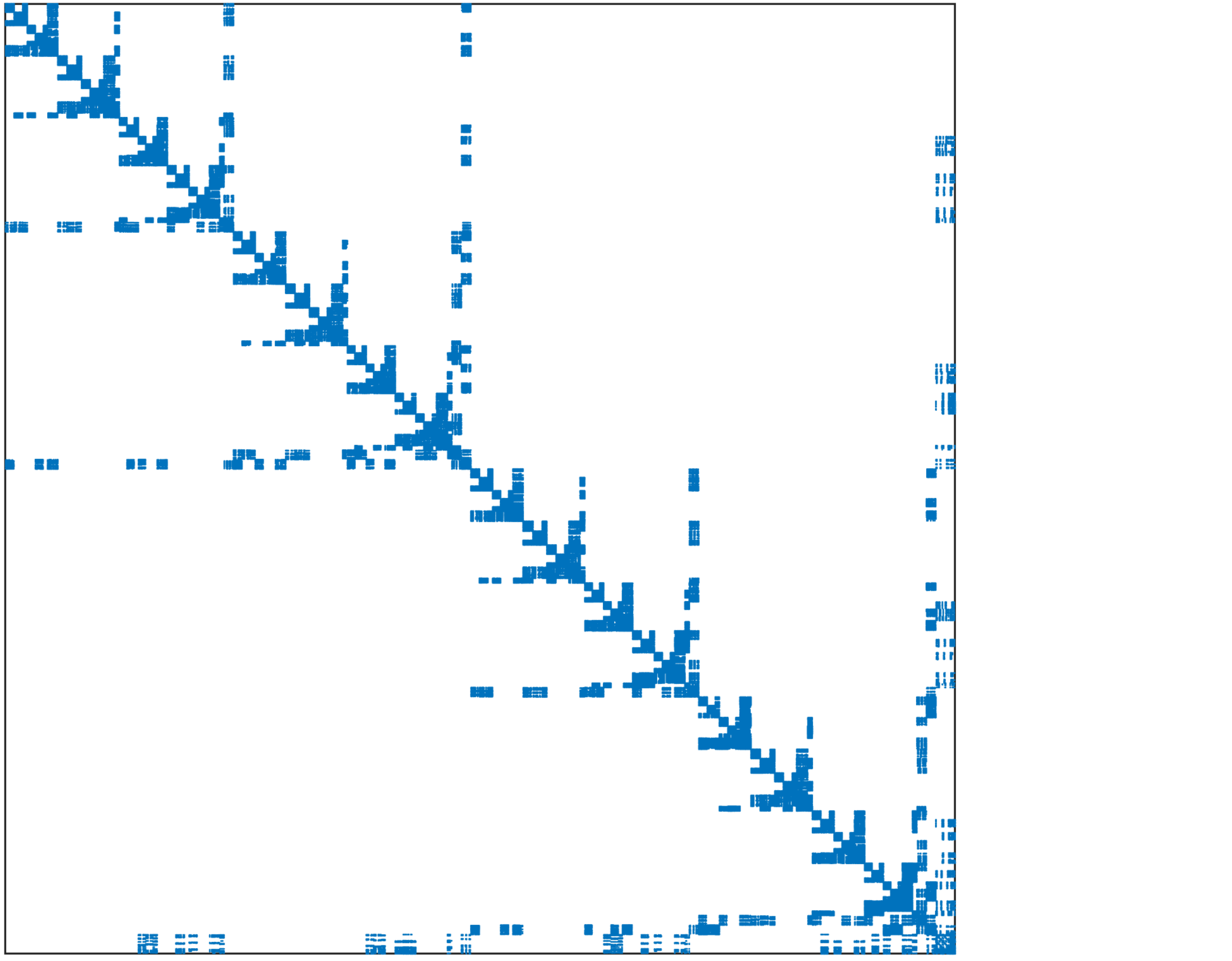}\end{overpic}
		\caption{rIGA reordered matrix: $\nnz=\,$1,217,968}
		\label{fig.MatPatterns.rIGAreorder}
	\end{subfigure}
	\begin{subfigure}{0.33\textwidth} 
		\begin{overpic}[width=1\linewidth,trim={0cm 0cm 7.25cm 0cm},clip] {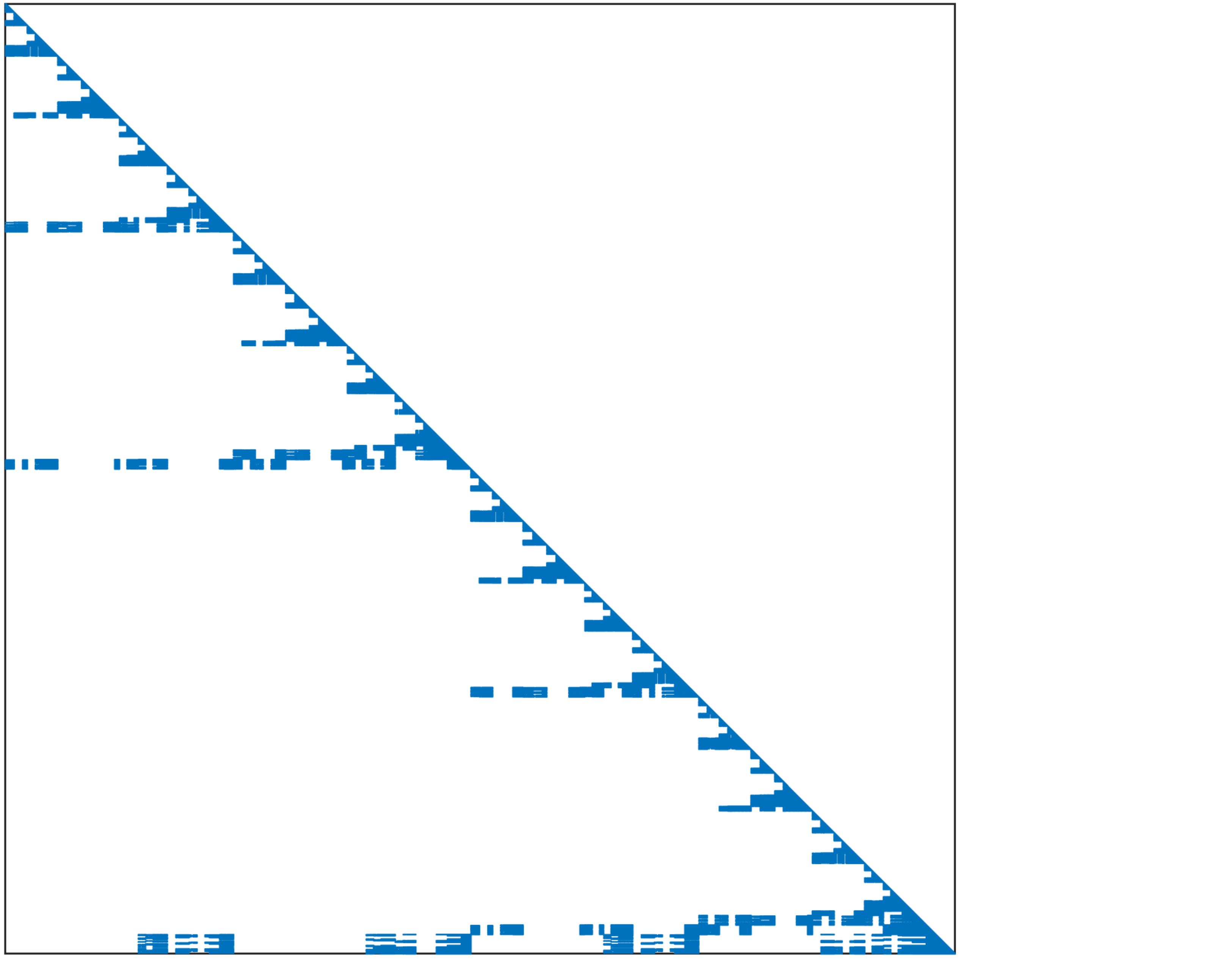}\end{overpic}
		\caption{rIGA lower factor: $\nnz=\,$2,117,582}
		\label{fig.MatPatterns.rIGAfactor}
	\end{subfigure}
  \caption{Matrix patterns for a ${64\times64}$ domain discretization using IGA and rIGA \ibH(\curl) spaces with ${\p=4\.}$: (a,d) original system matrices, (b,e) reordered matrices obtained by nested dissection permutation, and (c,f) lower triangular Cholesky factors. We obtain an rIGA discretization by performing two partitioning levels, creating ${16\times16}$ macroelements.}
	\label{fig.MatPatterns}
\end{figure}

We investigate how rIGA modifies the system matrices and the respective LU (or Cholesky) factors; thus, we consider an \ibH(\curl) space discretized by ${64\times64}$ elements and polynomial order ${\p=4}$. \figs{\ref{fig.MatPatterns.IGAmatrix}}{\ref{fig.MatPatterns.rIGAmatrix}} show the IGA and rIGA matrix patterns, respectively. We use two levels of partitioning that creates ${16\times16}$ macroelements. We use the nested dissection permutation~\cite{ Karypis1998} to obtain the reordered matrices (\figs{\ref{fig.MatPatterns.IGAreorder}}{\ref{fig.MatPatterns.rIGAreorder}}). An rIGA-discretized system has more number of degrees of freedom and, hence, a higher number of nonzeros ($\nnz$). However, the lower continuity of the mesh at hyperplanes leads to lower interconnections in the matrix connectivity graph. Thus, the factorization is obtained faster and the respective factors have a lower number of nonzero terms (\figs{\ref{fig.MatPatterns.IGAfactor}}{\ref{fig.MatPatterns.rIGAfactor}}\.).


\section{Numerical solution to quadratic eigenvalue problems} 
\label{sec.EigsolutionAlgorithms}

We review the numerical solution approach to quadratic eigenvalue problems. For the sake of brevity and to follow the same notation for both model problems introduced in~\Sec{\ref{sec.ModelProblems}}, we denote the stiffness, damping, and mass matrices by $\K$, $\C$, and $\M$, respectively. We also drop the subscript $h$ and consider ${(\lm,\tuu)}$ as the approximate eigenvalue and eigenvector pair.  We review the generalized eigenproblems since we first linearize the quadratic eigenvalue problem and then solve a generalized eigensystem (see, e.g.,~\cite{ Olson1989, Bermdez2000, Tisseur2001, Campos2016}).

\subsection{Generalized eigenproblems}
\label{sub:LinearEig}

Let us consider the following generalized eigenproblem, 
\begin{align}
(\A-\lm\B)\.\tuu=\textbf{0}\,,
\label{eq.LinearEig}
\end{align}
characterized by a linear {\em matrix pencil} ${\Lr(\lm)\ceq\A-\lm\B\.}$.  The term generalized distinguishes~\eqref{eq.LinearEig} from the standard eigenproblem for which $\B$ is the unit matrix.  If both $\A$ and $\B$ are Hermitian matrices and $\B$ is definite, $\Lr(\lm)$ is a Hermitian definite pencil, and all eigenvalues are real; otherwise, some eigenvalues may become complex.  When solving QEPs, via linearization, the linear pencil, on most occasions, is non-Hermitian or indefinite, resulting in complex eigenvalues.
For linear Hermitian definite pencils, most eigensolvers obtain a portion of the eigenspectrum by executing a {\em sequence} of numerical algorithms, among which are: shift-and-invert spectral transformation, Arnoldi (or Lanczos) decomposition, restarting, spectrum slicing, and deflation (see, e.g.,~\cite{ Hashemian2021}). However, not all algorithms apply to the quadratic case as we may deal with non-Hermitian (or indefinite) linearized pencils.  In the following, we present a short review of those algorithms essential to the efficiency study of the eigencomputation of QEPs.


\subsubsection{Shift-and-invert spectral transformation} 
\label{ssub:ShiftInvert}

When solving eigenproblem~\eqref{eq.LinearEig}, an efficient way of computing eigenpairs is to shift the spectrum by ${\sh\in\CC}$ and solve the shift-and-invert spectral transformed eigenproblem\footnote{We first transform~\eqref{eq.LinearEig} to ${(\A-\sh\.\B)\.\tuu=(\lm-\sh)\.\B\tuu}$ and, then, invert it to reach~\eqref{eq.ShiftInvert}.} (see, e.g.,~\cite{ Ericsson1980, NourOmid1987, Grimes1994, XUE2011, Demmel2000}),
\begin{align}
(\A-\sh\.\B)^{-1}\B\tuu=\theta\.\tuu\,,
\label{eq.ShiftInvert}
\end{align}
where ${\theta=1/(\lm-\sh)}$.  In large multidimensional eigensystems with multiple clustered eigenvalues, this approach is fast and accurate when calculating eigenvalues near the shift $\sh$.  By defining the {\em operator} matrix, 
\begin{align}
\S\ceq(\A-\sh\.\B)^{-1}\B=\Lr(\sh)^{-1}\B\,,
\label{eq.OperatorMatrix}
\end{align}
we express the shift-and-invert eigenproblem as ${\S\.\tuu=\theta\.\tuu}$.


\subsubsection{Arnoldi decomposition} 
\label{ssub:Arnoldi}

To solve the eigenproblem~\eqref{eq.ShiftInvert}, when $\B$ is positive definite, we employ the Arnoldi decomposition, which can deal with either Hermitian or non-Hermitian $\A$.  For a given shift $\sh$ and either Cholesky or LU factorization of $\Lr(\sh)$, the $m$-step Arnoldi decomposition,
\begin{align}
\S\.\V_m=\V_m\H_m+\beta\.\vv_{m+1}\e_m^T\,,
\label{eq.Arnoldi}
\end{align}
projects the eigensystem~\eqref{eq.ShiftInvert} of size $\N$ onto Krylov subspace of size $m$~${(m\ll N)}$ by reducing the ${\N\times \N}$ operator matrix~$\S$ to an ${m\times m}$ upper Hessenberg matrix $\H_m$. Herein, $\e_m$ is the ${m\.}$th coordinate vector and the term ${\beta\.\vv_{m+1}\e_m^T}$ is the residual of the $m$-step Arnoldi decomposition (see, e.g.,~\cite{ Demmel2000, Stewart2001}).
In the above equation,
\begin{align}
\H_m\ceq\begin{bmatrix}
h_{1,1} & h_{1,2} & \cdots & h_{1,m-2}   & h_{1,m-1}   & h_{1,m}  \\
h_{2,1} & h_{2,2} & \cdots & h_{2,m-2}   & h_{2,m-1}   & h_{2,m}  \\
        & h_{3,2} & \cdots & h_{3,m-2}   & h_{3,m-1}   & h_{3,m}  \\
        &         & \ddots & \vdots      & \vdots      & \vdots   \\
        &         &        & h_{m-1,m-2} & h_{m-1,m-1} & h_{m-1,m}\\
        &         &        &             & h_{m,m-1}   & h_{m,m}
\end{bmatrix},
\label{eq.ArnoldiTm}
\end{align}
and ${\V_m\ceq[\vv_1,\vv_2,...,\vv_m]}$ is the ${\N\times m}$ matrix of Arnoldi basis vectors.
We obtain the components of $\H_m$ and Arnoldi bases using the recurrence given by \Alg{\ref{alg.Arnoldi}}.

\begin{algorithm}[!h]
\caption{Arnoldi decomposition for solving (linear) generalized eigenproblems.}
\label{alg.Arnoldi}
\textbf{Input:} Positive definite matrix $\B$, factorized form of $\Lr(\sh)$, size of Krylov subspace $m$, and initial vector $\vv_1$\,.\\
\textbf{Output:} Arnoldi bases $\vv_j$ $(j=1,2,...,m+1)$\., upper Hessenberg matrix $\H_m$, and the residual factor $\beta$\,.
\begin{spacing}{1.15}
\begin{algorithmic}[1]
\State Normalize: $\vv_1\leftarrow\vv_1/\normB{\vv_1}$\,,\quad $\V_1=\vv_1$
\ForAll{$j=1,2,...,m$,}
\State $\rr=\S\vv_j=\Lr(\sh)^{-1}\B\vv_j$					\label{alg.Arnoldi.rcalc}
\State $\hh=\V_j^T\B\rr$\,,\quad $h_{1:j,j}\leftarrow\hh$ 	\label{alg.Arnoldi.hcalc}
\State $\zz=\rr-\V_j\.\hh$ 									\label{alg.Arnoldi.GSchmidt}
\State $h_{j+1,j}=\normB{\.\zz\.}$
\State $\vv_{j+1}=\zz/h_{j+1,j}$\,,\quad $\V_{j+1}\leftarrow\begin{bmatrix}\V_j &\vv_{j+1}\end{bmatrix}$
\EndFor
\State $\beta=h_{m+1,m}$
\end{algorithmic}
\end{spacing}
\end{algorithm}

The Arnoldi vector $\vv_{m+1}$ is $\B$-orthogonal with respect to the columns of $\V_m$ in the Gram--Schmidt sense, resulting in ${\V_m^T\.\B\vv_{m+1}=0}$.  Hence, the $\B$-inner product of $\V_m$ premultiplied in~\eqref{eq.Arnoldi} leads to the following equation, noting that~$\V_m$ is $\B$-orthogonal (i.e., ${\V_m^T\.\B\V_m=\I}$\.)\.:
\begin{align}
\V_m^T\.\B\.\S\V_m=\H_m\,.
\label{eq.ArnoldiRitz}
\end{align}
In~\eqref{eq.ArnoldiRitz}, $\H_m$ is the $\B$-orthogonal projection of $\S$ onto Krylov subspace. Therefore, the~$m$ eigenpairs of~$\H_m$ (referred to as Ritz pairs) are the Rayleigh--Ritz approximation of the $m$ eigenpairs ${(\theta_j,\tuu_j)}$ of the shift-and-invert problem~\eqref{eq.ShiftInvert}.
Since $\H_m$ is not symmetric, some eigenvalues of~\eqref{eq.ShiftInvert} may be complex. 

\begin{remark}
  If $\A$ is Hermitian and ${\sh\in\R}$, the calculation of $\hh$ in line~\ref{alg.Arnoldi.hcalc} of~\Alg{\ref{alg.Arnoldi}} includes only the last basis vector (i.e.,~${h_{j,j}=\vv_j^T\B\rr}$) and the Gram--Schmidt orthogonalization (line~\ref{alg.Arnoldi.GSchmidt}) is only performed against the last two basis vectors (i.e., ${\zz=\rr-h_{j-1,j}\vv_{j-1}-h_{j,j}\vv_j}$). In this case, $\H_m$ converts to a symmetric tridiagonal matrix resulting in real eigenvalues. \Eq{\eqref{eq.Arnoldi}} is then referred to as the Lanczos decomposition.
\label{rem.Lanczos}
\end{remark}


\subsubsection{Restarting} 
\label{ssub:Restarting}

When seeking a large number of eigenvalues ($\Nev$), the norm of the residual term in~\eqref{eq.Arnoldi} grows, and the approximation quality of Ritz pairs deteriorates.  Thus, Krylov eigensolvers incorporate iterative restarting mechanisms to maintain the residual norm lower than the desired tolerance $\epsilon$.  By setting $m$ slightly larger than $\Nev$\., the eigensolver keeps a set of {\em good} eigenpairs after each Arnoldi process ---\.those with residual norms below $\epsilon$. Then, a new $m$-step Arnoldi process restarts in the next iteration, benefiting from the previously obtained spectral approximation.  The restarting process continues until we compute all requested eigenpairs.  Well-known restarting techniques are the thick-restart Lanczos method~\cite{ Wu2000}, and its nonsymmetric equivalent for the Arnoldi case referred to as the Krylov--Schur method~\cite{ Stewart2002, Stewart2002a}.  We use the latter in our eigencomputations, but we omit details for the sake of brevity.


\subsubsection{Spectrum slicing} 
\label{ssub:SpectrumSlicing}

The spectrum slicing technique is only applicable when all eigenvalues are real~\cite{ Grimes1994, Campos2012}. Thus, most quadratic eigensolvers cannot use this technique as their spectrum may be complex\footnote{Under specific circumstances described in~\cite{ Campos2020} the spectrum slicing technique applies to QEPs. Still, we do not study these cases here.}. Nevertheless, we briefly review this algorithm to highlight a significant difference between generalized and quadratic eigenproblems.

By definition, if all eigenvalues of an arbitrary matrix $\A$ are real, the {\em inertia} of $\A$ is given by three integers ${\nu\.(\A)}$, $\zeta(\A)$, and ${\pi\.(\A)}$ denoting the numbers of negative, zero, and positive eigenvalues of $\A$, respectively.  For linear Hermitian definite pencils, we consider an LDL Cholesky factorization of $\Lr(\sh)$, which is 
${\A-\sh\.\B=\LL\DD\LL^T},$
being $\DD$ a~diagonal matrix.  Then, based on Sylvester's law of inertia~\cite{ Parlett1998}, one writes ${\nu\.(\A-\sh\.\B)=\nu\.(\DD)}$, implying that the number of eigenvalues of $\Lr(\lm)$ smaller than $\sh$ is equal to the number of negative eigenvalues of $\DD$.  As a result, considering two shifts $\sh_k$ and $\sh_{k+1}$, the spectrum subinterval $[\sh_k,\sh_{k+1}]$ has ${\nu\.(\A-\sh_{k+1}\B)-\nu\.(\A-\sh_k\B)}$ eigenvalues including their multiplicities.  The inertia information helps to {\em slice} the spectrum efficiently into subintervals, leading to the efficient spectrum slicing technique.  When $\Nev$ is large, spectrum slicing has two main benefits. First, the eigensolution accuracy increases because we look for eigenvalues in a smaller spectrum subinterval (i.e., eigenvalues are closer to the shift). Second, spectrum slicing prevents the deterioration of the convergence rate of the Arnoldi decomposition~\eqref{eq.Arnoldi} because we project onto a smaller Krylov subspace at each Arnoldi run.  In this manner, eigensolvers set the subspace size $m$ to be slightly higher than the number of to-be-computed eigenvalues per shift (see, e.g.,~\cite{ Hashemian2021}).  Each shift then contains the corresponding Arnoldi recurrences and restarts described in previous sections.
Importantly, when dealing with the non-Hermitian or indefinite pencils that result from the linearization of quadratic eigenproblems, the lack of exact inertia (due to a complex spectrum) prevents eigensolvers to perform this technique~\cite{ Nakatsukasa2019}. For such eigensystems, we select only one {\em target} shift inside the spectrum region of interest.


\subsection{Quadratic eigenproblems} 
\label{sub:QuadEig}

We consider the quadratic eigenproblem,
\begin{align}
(\K+\lm\C+\lm^2\M)\.\tuu=\textbf{0} \,,
\label{eq.QuadEig}
\end{align}
characterized by the quadratic matrix pencil ${\Qr(\lm)\ceq\K+\lm\C+\lm^2\M}$\..  For ${\N\times\N}$ system matrices, this eigenproblem has $2\N$ eigenpairs.  In many practical cases, $\Qr(\lm)$ is either Hermitian or gyroscopic (c.f., our model problems in~\Sec{\ref{sec.ModelProblems}}), which implies that the eigenvalues may appear in a combination of complex conjugate pairs and real values.  Under special circumstances that all eigenvalues of a quadratic Hermitian eigensystem may be purely real, we call $\Qr(\lm)$ a {\em hyperbolic} pencil~\cite{ Guo2005, Niendorf2010}. But this particular case is not studied here.


\subsubsection{Linearization} 
\label{ssub:Linearization}

Eigensolvers perform different types of projection for computing the eigenvalues of quadratic eigenproblems.  Methods such as Jacobi--Davidson~\cite{ Sleijpen1996} and second-order Arnoldi~(SOAR)~\cite{ Bai2005} directly project the quadratic eigenproblem.  However, a common way to solve~\eqref{eq.QuadEig} is to perform a linearization by projecting an equivalent linearized eigenproblem of a doubled size onto Krylov subspace (see, e.g.,~\cite{ Olson1989, Bermdez2000, Tisseur2001, Campos2016}).  Some linearizations preserve the Hermitian structure of matrices~\cite{ Campos2020}. Nevertheless, the obtained linear pencil is indefinite in most cases.  Herein, we use the following linearization of~\eqref{eq.QuadEig}, which applies to our two model problems:
\begin{align}
\Bigg(\underbrace{\begin{bmatrix}\textbf{0} & \I\\ -\K & -\C \end{bmatrix}}_{\A}
  -\lm\underbrace{\begin{bmatrix}\I & \textbf{0}\\ \textbf{0} &\M\end{bmatrix}}_{\B}\Bigg)
   \begin{bmatrix}\tuu\\ \lm\tuu\end{bmatrix}=\textbf{0}\,,
\label{eq.Linearization}
\end{align}
where $\A$ and $\B$ are doubled-size matrices. After linearizing the quadratic eigenproblems, we use the same numerical algorithms as the generalized eigensystems. In the resulting linearized eigenproblem~\eqref{eq.Linearization}, matrix $\B$ is positive definite, while matrix $\A$ is non-Hermitian. Thus, the Arnoldi recurrence of \Alg{\ref{alg.Arnoldi}} holds. However, we deal with $2\N\times2\N$ matrices that will increase the computational cost of the eigensolution.  We may resolve this issue by exploiting the block structure of the linearization~\eqref{eq.Linearization} without explicitly constructing the doubled-size matrices $\A$ and $\B$ (see, e.g.,~\cite{ Lu2016, Campos2020}). Considering
\begin{align}
\Lr(\sh)=\A-\sh\.\B=
\begin{bmatrix}-\sh\.\I & \I \\ -\K & -\C-\sh\.\M\end{bmatrix}\,,
\label{eq.BlockLinearized}
\end{align}
using the block factorization of $\Lr(\sh)$ and taking into account the definition of the quadratic pencil $\Qr(\sh)$, we write the operator matrix as
\begin{align}
\S=\Lr(\sh)^{-1}\B=
\begin{bmatrix}\Qr(\sh)^{-1} & \textbf{0} \\\textbf{0}& \Qr(\sh)^{-1} \end{bmatrix}
\begin{bmatrix}-\C-\sh\.\M   	 & -\M        \\\K        & -\sh\.\M        \end{bmatrix}\,.
\label{eq.BlockOperator}
\end{align}
We then multiply the operator matrix with the Arnoldi basis vectors $\vv$ (line~\ref{alg.Arnoldi.rcalc} of \Alg{\ref{alg.Arnoldi}}) by splitting $\vv$ into two upper and lower $\N\times1$ blocks as ${\vv=\smaller{\begin{bmatrix}\vv_u\\\vv_l\end{bmatrix}}}$\.. Hence, by writing ${\rr=\smaller{\begin{bmatrix}\rr_u\\\rr_l\end{bmatrix}}}$ in a similar way, we obtain
\begin{align}
\rr=\Lr(\sh)^{-1}\B\vv=\begin{bmatrix}
-\Qr(\sh)^{-1}(\sh\.\M\vv_u+\C\vv_u+\M\vv_l) \\\sh\.\rr_u+\vv_u\end{bmatrix}\..
\label{eq.BlockOperatorVector}
\end{align}
In this context, we write every $\B$-norm calculation as:
\begin{align}
\normB{\vv}^2=\vv^T\B\vv =
\begin{bmatrix}\vv_u^T & \vv_l^T\end{bmatrix}
\begin{bmatrix}\I & \textbf{0}\\\textbf{0} & \M\end{bmatrix}
\begin{bmatrix}\vv_u \\ \vv_l\end{bmatrix} =
\vv_u^T\vv_u+\vv_l^T\M\vv_l
\..
\end{align}

\begin{algorithm}[!pb]
	\caption{Numerical eigenanalysis procedure for solving quadratic eigenproblems.}
	\label{alg.QuadEig}
	\textbf{Input:} ${\N\times\N}$ system matrices $\K$, $\C$ and $\M$, target shift $\sh$, requested number of eigenvalues $\Nev$\., size of Krylov subspace ${m>\Nev}$\., initial vector $\vv_1$\., and residual tolerance $\epsilon$.\\
	\textbf{Output:} eigenpairs $\lm$ and $\tuu$\,.
	\begin{spacing}{1.15}
		\begin{algorithmic}[1]
			\State Normalize: $\vv_1\leftarrow\vv_1/\normB{\vv_1}$\,,\quad $\V_1=\vv_1$
			\State Factorize $\Qr(\sh)$										\label{alg.QuadEig.fact}
			\State $\Nit=0$\, \quad /* a counter for total number of iterations */
			\State ${R=10^{\.6}}$ \quad /* setting a big value for residual norm to get into iterations */
			\While{$R>\epsilon$,}
			\State $\Nit\leftarrow\Nit+1$
			\ForAll{$j=1,2,...,m$,} \quad /* $m$-step Arnoldi recurrence */ \vspace{2pt}
			\State $\rr=\smaller{\begin{bmatrix}
				-\Qr(\sh)^{-1}\big(\sh\.\M\vv_{j,u}+\C\vv_{j,u}+\M\vv_{j,l}\big) \\
				\sh\.\rr_u+\vv_{j,u}\end{bmatrix}}$							\label{alg.QuadEig.fb}
			\State $\hh=\V_j^T\B\rr=\smaller{\begin{bmatrix}
				\vv_{1,u}^T\rr_u+\vv_{1,l}^T\.\M\rr_l\\[-3pt]\vdots\\
				\vv_{j,u}^T\rr_u+\vv_{j,l}^T\.\M\rr_l
				\end{bmatrix}}$\,, \quad $h_{1:j,j}\leftarrow\hh$ 			\label{alg.QuadEig.mv}
			\State $\zz=\rr-\V_j\.\hh$ 										\label{alg.QuadEig.Vv}
			\State $h_{j+1,j}=\normB{\.\zz\.}=
			\big(\zz_u^T\zz_u+\zz_l^T\M\zz_l\big)^{1/2}$				\label{alg.QuadEig.vv}
			\State $\vv_{j+1}=\zz/h_{j+1,j}$\,,\quad $\V_{j+1}\leftarrow\begin{bmatrix}\V_j &\vv_{j+1}\end{bmatrix}$
			\EndFor
			\State $\beta=h_{m+1,m}$
			\State Construct matrix $\H_m$ and compute Ritz pairs
			\State Compute numerical eigenpairs $\lm_j$ and $\tuu_j$ $(j=1,2,...,m)$
			\State Compute largest residual norm $R$
			\State Keep adequate numbers of eigenpairs and compute initial values for restart
			\EndWhile
		\end{algorithmic}
	\end{spacing}
\end{algorithm}

\noindent
This block representation allows us to perform factorization and matrix--vector operations of the Arnoldi recurrence using ${\N\times\N}$ system matrices. Other steps of solving quadratic eigenproblems exploit the same block-product structure. \Alg{\ref{alg.QuadEig}} presents a pseudo-code describing the numerical eigenanalysis procedure for solving quadratic eigenproblems. However, the entire process might be more complicated. For instance, the size of Krylov subspace and the adequate number of Arnoldi bases to keep before restart play essential roles in the efficiency of the eigensolution. These variables mainly depend on the size of the eigensystem, the total number of requested eigenpairs, and the selected eigensolver package.


\section{Eigencomputation cost} 
\label{sec.CostEig}

\subsection{General overview} 
\label{sub:CostOverview}

The total eigencomputation cost of quadratic eigenproblems is mostly governed by the following operation sets
\begin{itemize}
\item LU factorization (\fa) of ${\Qr(\sh)}$;
\item Forward/backward (\fb) eliminations (i.e., multiplications of LU factors by vectors);
\item Matrix--vector (\mv) and vector--vector (\Vv) products in the sense of Krylov projection (see~\Alg{\ref{alg.QuadEig}}).
\end{itemize}
We consider other costs (e.g., those related to creating system matrices) as lower-order terms and exclude them from our cost estimation. Thus, we write
\begin{align}
\Cr \approx \Nfa\.\Cr_\fa + \Nfb\.\Cr_\fb + \Nmv\.\Cr_\mv + \Nvv\.\Cr_\Vv\,,
\label{eq.Cost}
\end{align}
where with $\Cr$ and $\Nr$, we refer to the cost and number of operations, respectively. 

For linear Hermitian definite pencils, the spectrum slicing entails performing multiple factorizations (one per shift). At the same time, thanks to a smaller Krylov subspace, fewer matrix--vector operations are required (see \Sec{\ref{ssub:SpectrumSlicing}}).  Additionally, in such eigensystems, only the last two basis vectors participate in each recurrence step (see \Rem{\ref{rem.Lanczos}}), implying that the cost of vector--vector products is negligible.  As a result, for moderate to large generalized eigensystems (characterized by linear Hermitian definite pencils), matrix factorization governs the total eigencomputation cost (see~\cite{ Hashemian2021}).  On the other hand, for quadratic eigenproblems with non-Hermitian (or indefinite) linearized pencils, we perform one shift and, therefore, one factorization for all requested eigenpairs, resulting in a larger Krylov subspace (the subspace size grows as $\Nev$ increases). Moreover, a full basis matrix~$\V$ participates in the Arnoldi recurrence.  Thus, when solving moderately sized problems, the other three operations in~\eqref{eq.Cost} govern the eigencomputation cost. Their contribution to the total cost is remarkable when we look for a large number of eigenpairs. Nevertheless, the matrix factorization is asymptotically the most expensive operation as the problem size grows.


\subsection{Eigencomputation cost per operation} 
\label{sub:CostOperations}

\subsubsection{Matrix factorization} 
\label{ssub:CostFact}

As we state above, solving quadratic eigenproblems entails one LU factorization of the quadratic pencil ${\Qr(\sh)}$ for the target shift $s$ inside the eigenspectrum region of interest (i.e., ${\Nfa=1}$, see line~\ref{alg.QuadEig.fact} of \Alg{\ref{alg.QuadEig}}). 
To determine the factorization cost, we use the theoretical estimates of IGA and rIGA discretizations for \ibH(\curl) and \ibH(\div) spaces~\cite{ Garcia2019}. Thus, for 2D problems, we write
\begin{align}
\Cr_{\fa\.,\.\IGA} 	&\approx \Or\big(\N^{1.5}(\p-1)^3\big)\,, 	\label{eq.CostFactIGA} \\
\Cr_{\fa\.,\.\rIGA} &\approx \Or\big(\N^{1.5}(\p-1)\big)\,. 	\label{eq.CostFactrIGA}
\end{align}
The number of degrees of freedom of an rIGA system is slightly higher than its IGA counterpart (a.k.a., ${\N_\rIGA\approx\N_\IGA}$). Thus, we consider ${\N=\N_\IGA}$ in both~\eqref{eq.CostFactIGA}~and~\eqref{eq.CostFactrIGA} and obtain a cost improvement of ${\Or\big((\p-1)^2\big)}$ for large problem sizes when using rIGA discretization with optimal macroelements containing 16 elements in each direction. For smaller eigensystems, computational savings above ${\Or(\p-1)}$ are hardly obtainable (see~\cite{Garcia2019}).


\subsubsection{Forward/backward elimination} 
\label{ssub:CostFb}

The number of forward/backward eliminations is equal to the size of Krylov subspace, $m$, multiplied by the number of iterations, $\Nit$ (i.e., ${\Nfb=m\.\Nit}$\., see line~\ref{alg.QuadEig.fb} of \Alg{\ref{alg.QuadEig}} in computing the upper block of~$\rr$). 
The forward/backward elimination cost is proportional to the number of nonzero terms of the LU factors, which we estimate for 2D IGA and rIGA discretizations as follows (see~\cite{Garcia2017,Hashemian2021})\.:
\begin{align}
\Cr_{\fb\.,\.\IGA}  &\approx \Or\big(\N(\p-1)^2\big)\,, 	\label{eq.CostFBIGA} \\
\Cr_{\fb\.,\.\rIGA} &\approx \Or\big(\N(\p-1)\big)\,. 		\label{eq.CostFBrIGA}
\end{align}
Again, the cost reduction of ${\Or(\p-1)}$ with respect to IGA in~\eqref{eq.CostFBrIGA} is asymptotically apparent only for large problem sizes and when employing macroelements with 16 elements in each direction.


\subsubsection{Matrix--vector multiplication} 
\label{ssub:CostMv}

The quadratic eigenproblem solver calls the matrix--vector operator in lines~\ref{alg.QuadEig.fb},~\ref{alg.QuadEig.mv}~and~\ref{alg.QuadEig.vv} of \Alg{\ref{alg.QuadEig}}. Thus, one obtains ${\Nmv=5m\.\Nit}$ for the participation of either $\M$ or $\C$ in the process. This numerical operation has a cost proportional to the number of nonzero entries, $\nnz$\., of system matrices ---\.that is related to the sum of interactions of each basis function with all other bases~\cite{Collier2013}. For 2D eigensystems, we write
\begin{align}
\Cr_{\mv\.,\.\IGA}  &= \Or\big(\nnz(\M_\IGA)\big)\approx\Or(\N_\IGA\p^2)\,, 	\label{eq.CostMVIGA} \\
\Cr_{\mv\.,\.\rIGA} &= \Or\big(\nnz(\M_\rIGA)\big)\approx\Or(\N_\rIGA\p^2)\,. 	\label{eq.CostMVrIGA}
\end{align}
An rIGA-discretized system has a slightly higher number of degrees of freedom. Thus, matrix--vector products suffer a slight degradation when using rIGA. We distinguish between $\N_\IGA$ and $\N_\rIGA$ in~\eqref{eq.CostMVIGA}~and~\eqref{eq.CostMVrIGA} since it is the only difference in the cost of this numerical operation with respect to the employed discretization.


\subsubsection{Vector--vector multiplication} 
\label{ssub:CostVv}

The vector--vector product is theoretically the cheapest numerical operation of eigenanalysis. It has a computational complexity proportional to the number of entries of vectors, $\N$, which is significantly smaller than $\nnz(\M)$. However, when solving quadratic eigenproblems, we perform a large number of vector--vector products in the sense of Gram--Schmidt orthogonalization (lines~\ref{alg.QuadEig.fb}\.--\ref{alg.QuadEig.vv} of \Alg{\ref{alg.QuadEig}}). In particular, at the ${j\.}$-th Arnoldi step~${(j=1,2,...,m)}$, the participation of $\V_j$ incorporates with $2j$ multiplications of vectors of size $\N$ in line~\ref{alg.QuadEig.mv}, and~$2\N$ products of vectors of size $j$ in line~\ref{alg.QuadEig.Vv}. 
At each iteration, the eigensolver repeats these lines $m$ times during the Arnoldi process. Thus, one obtains ${\Nvv\approx2m^2\Nit}$ as the total number of multiplications of vectors of size $\N$. As a result, the vector--vector products have a comparable cost as matrix--vector operations when we seek to compute a large number of eigenvalues (i.e., a large~$m$). 
In contrast to other operations, the cost of each vector--vector product is not directly related to the polynomial degree~$\p$ as it only depends on $\N$\.:
\begin{align}
\Cr_{\Vv\.,\.\IGA}  &= \Or(\N_\IGA)\,, 		\label{eq.CostVVIGA} \\
\Cr_{\Vv\.,\.\rIGA} &= \Or(\N_\rIGA)\,. 	\label{eq.CostVVrIGA}
\end{align}


\subsection{Summary of eigencomputation cost} 
\label{sub:CostSummary}

\Eqs{\eqref{eq.CostFactIGA}}{\eqref{eq.CostFactrIGA}} show that matrix factorization is the most expensive operation when $\N$ is large. The other three operations in~\eqref{eq.Cost} are not inherently as expensive as the matrix factorization. However, the eigensolver calls them many times during the Krylov projection. Thus, they become decisive when the number of requested eigenvalues $\Nev$ is large.  Since the subspace size~$m$ is slightly higher than $\Nev$\., we assume ${m\approx\Nev}$ and summarize the contribution of each numerical operation to the total eigencomputation cost as tabulated in \tab{\ref{tab:CostSummary}}.

\begin{table}[!h]
	\centering
	\caption{Most expensive operations of the quadratic eigensolver and their contribution to the total eigencomputation cost for 2D eigensystems. Cost improvements of rIGA are under the assumption of being $\N$ sufficiently large and using optimal macroelements with 16 elements in each direction. Unless stated otherwise, ${\N=\N_\IGA}$\..}
	\label{tab:CostSummary}
	\small
	\begin{tabular}{@{}llllll@{}}
		\toprule
		Numerical operation                    	& 		& LU factorization    & FB elimination      & Mat--vec product          & Vec--vec product       \\[1.5pt]\toprule
		Number of times the operation is called &      	& 1                   & $\Nev\.\Nit$ 		& $5\Nev\.\Nit$  			& $2\Nev^2\.\Nit$ 		 \\\midrule
		Cost of performing one operation       	& IGA  	& $\Or(\N^{1.5}\p^3)$ & $\Or(\N\p^2)$       & $\Or(\N_\IGA\p^2)$      	& $\Or(\N_\IGA)$         \\[1pt]
												& rIGA 	& $\Or(\N^{1.5}\p)$   & $\Or(\N\p)$         & $\Or(\N_\rIGA\p^2)$     	& $\Or(\N_\rIGA)$        \\\midrule
		Total cost of performing the operation 	& IGA  	& $\Or(\N^{1.5}\p^3)$ & $\Or(\Nev\N\p^2)$   & $\Or(\Nev\N_\IGA\p^2)$  	& $\Or(\Nev^2\N_\IGA)$   \\[1pt]
												& rIGA 	& $\Or(\N^{1.5}\p)$   & $\Or(\Nev\N\p)$		& $\Or(\Nev\N_\rIGA\p^2)$ 	& $\Or(\Nev^2\N_\rIGA)$  \\\midrule
		Improvement/degradation of  			&      	& Improved by         & Improved by         & Degraded by 	            & Degraded by            \\ 
		performing the operation in rIGA		& 		& $\Or(p^2)$		  & $\Or(p)$			& $\N_\rIGA/\N_\IGA$		& $\N_\rIGA/\N_\IGA$	 \\\bottomrule
	\end{tabular}
\end{table}


\section{Numerical results} 
\label{sec.Results}

We report numerical results of our eigencomputations when using maximum-continuity and refined isogeometric analyses to solve the electromagnetic and vibroacoustic eigenproblems described in \Sec{\ref{sec.ModelProblems}}. 
We study the numerical efficiency and accuracy versus cost for different IGA and rIGA discretizations.
Before proceeding with results, we provide implementation details in the following subsection.


\subsection{Implementation details} 
\label{sub:Implementation}

We discretize our model problems using PetIGA-MF~\cite{petigamf}, which is a multifield extension of PetIGA~\cite{petiga} ---\.a high-performance isogeometric analysis framework based on PETSc~\cite{petsc}.
PetIGA-MF, utilized in many scientific and engineering applications (see, e.g.,~\cite{ Vignal2015, Espath2016, Cortes2017, Espath2017, Sarmiento2018, Garcia2019, clavijo2019, Hashemian20212}), uses different spaces for each field of interest and employs data management libraries to condense the data of multiple fields in a single object, thus simplifying the discretization construction. 

We also use SLEPc, the scalable library for eigenvalue problem computations~\cite{ slepc}, for performing the eigenanalysis.  SLEPc, implemented to solve different eigenproblems types (see, e.g.,~\cite{ Campos2012, Romero2014, Campos2016, Faber2018, Keeli2018, AraujoC2020}), allows us to employ the quadratic eigenproblem solver incorporating the shift-and-invert spectral transformation, linearization, the Arnoldi decomposition, and the Krylov--Schur restarting methods.

We use multifrontal direct solver MUMPS~\cite{ mumps} to perform matrix factorization and the forward/backward eliminations.  We employ the sequential version of MUMPS, which runs on a single thread (core).  We also use the automatic matrix reordering algorithm provided by METIS~\cite{ Karypis1998}.  


\subsection{Case studies}
\label{sub:CaseStudies}

\subsubsection{Electromagnetic wave propagation in a three-layer Earth model}
\label{ssub:EarthModel}

\begin{figure}[!pb]
	\centering
	\includegraphics{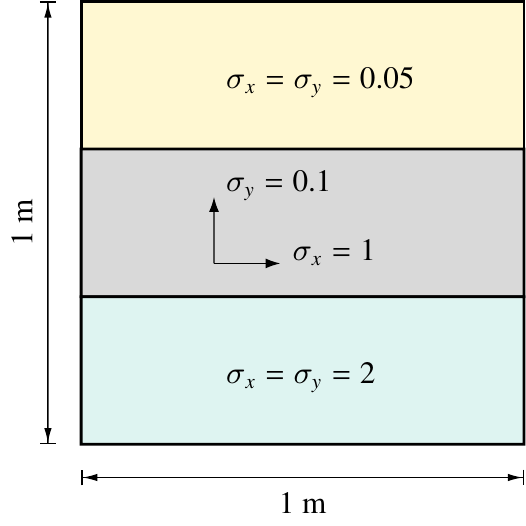}
	\caption{A three-layer Earth model as the physical domain of the electromagnetic eigenproblem. Each layer has electrical conductivities (given in ${\rm S/m}$) along horizontal and vertical directions.}
	\label{fig.EarthModel}
\end{figure}

We investigate the electromagnetic wave propagation eigenproblem. In particular, we consider a three-layer heterogeneous Earth model as a common case study in geosteering applications (see, e.g.,~\cite{ Hashemian20212, Shahriari2021}).  This 2D Earth model is essential in studying the full 3D wave propagation eigenproblem in transversely isotropic media~\cite{ Nicolet2006}.  As~\fig{\ref{fig.EarthModel}} shows, in each layer, we consider electrical conductivities along horizontal and vertical directions given by the conductivity matrix $\pmb{\sg}$ as follows:
\begin{align}
\pmb{\sg}\ceq\begin{bmatrix}\sg_x & 0 \\ 0 & \sg_y\end{bmatrix}\..
\label{eq.ConductivityMatrix}
\end{align}

\begin{figure}[!h]
	\centering
	\includegraphics{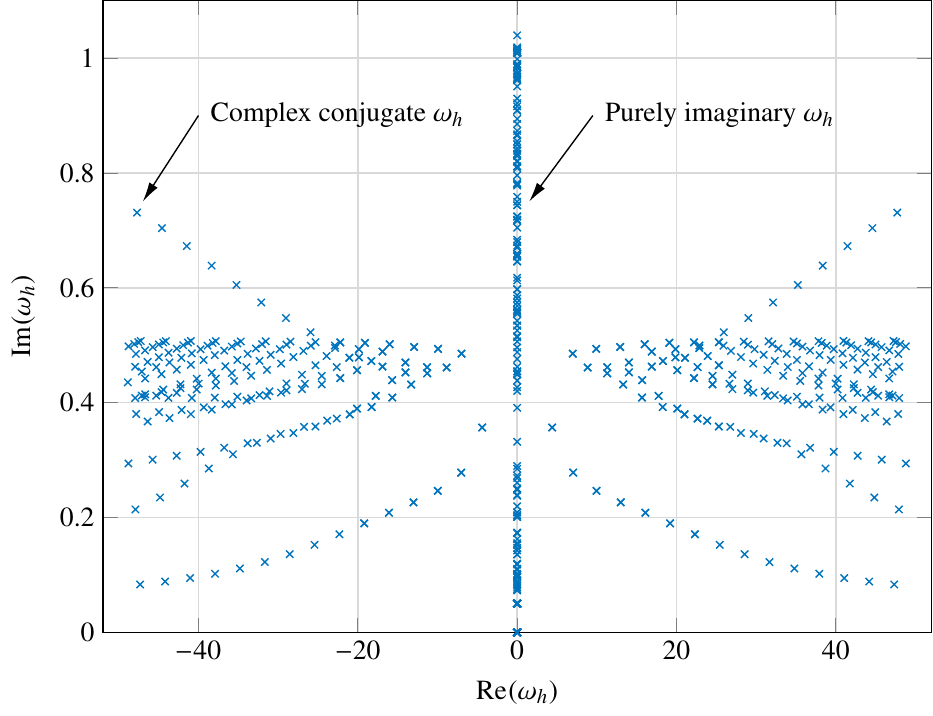}\vspace{-5pt}
	\caption{Approximate eigenvalues for electromagnetic wave propagation eigenproblem of the three-layer Earth model; $\omh$ values are given in ${\rm rad/s}$.}
	\label{fig.EarthModelEigs}
\end{figure}

\begin{figure}[!h]
	\begin{subfigure}{0.24\textwidth}
		\includegraphics{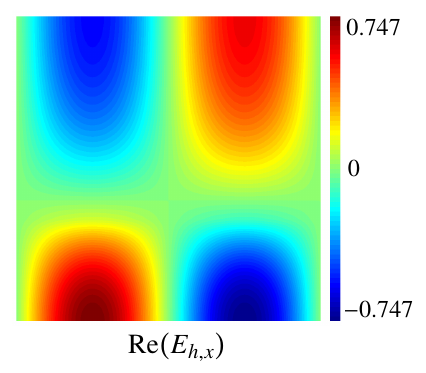}
	\end{subfigure}\hspace{3pt}
	\begin{subfigure}{0.24\textwidth}\centering
		\includegraphics{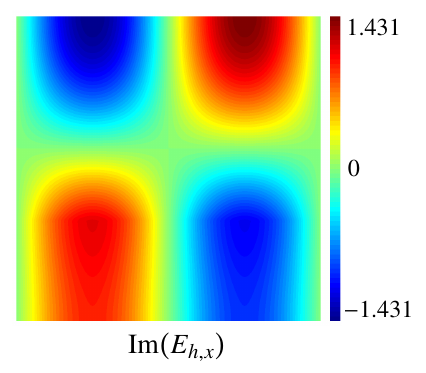}
	\end{subfigure}\hspace{3pt}
	\begin{subfigure}{0.24\textwidth}\centering
		\includegraphics{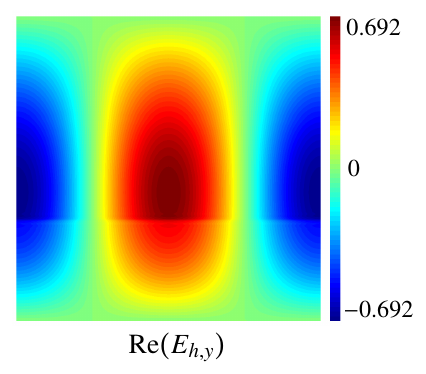}
	\end{subfigure}\hspace{3pt}		
	\begin{subfigure}{0.24\textwidth}\centering
		\includegraphics{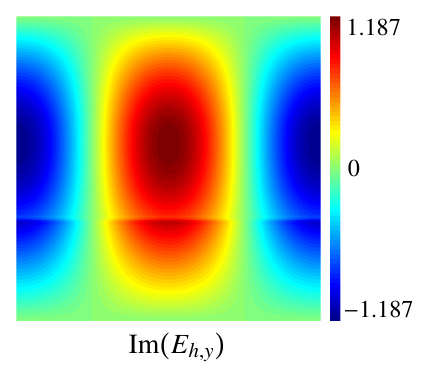}
	\end{subfigure}\\[-8pt]
	\begin{center}
		{\footnotesize (a) Approximate eigenfunctions associated with $\omh=7.024+0.282\.\ii~{\rm rad/s}$}
	\end{center} 
	\begin{subfigure}{0.24\textwidth}
		\includegraphics{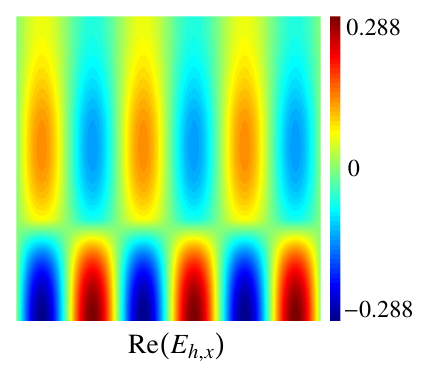}
	\end{subfigure}\hspace{3pt}
	\begin{subfigure}{0.24\textwidth}\centering
		\includegraphics{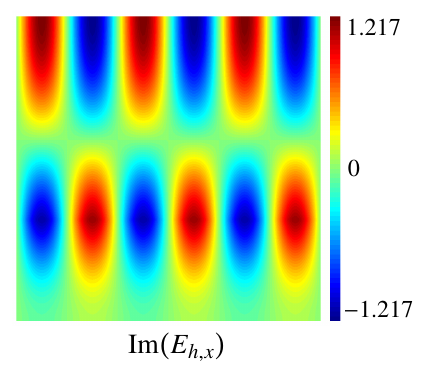}
	\end{subfigure}\hspace{3pt}
	\begin{subfigure}{0.24\textwidth}\centering
		\includegraphics{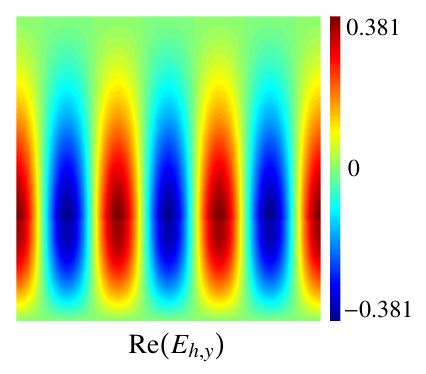}
	\end{subfigure}\hspace{3pt}		
	\begin{subfigure}{0.24\textwidth}\centering
		\includegraphics{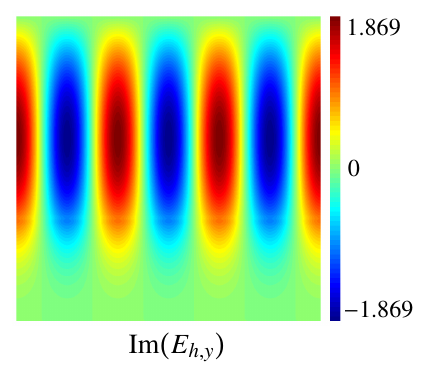}
	\end{subfigure}\\[-8pt]
	\begin{center}
		{\footnotesize (b) Approximate eigenfunctions associated with $\omh=19.200+0.193\.\ii~{\rm rad/s}$}
	\end{center}\vspace{-5pt}
	\caption{Two sample approximate eigenfunctions of the electromagnetic eigenproblem of~\fig{\ref{fig.EarthModel}}; $\EE_h$ values are given in ${\rm V/m}$.}
	\label{fig.EarthModelEigFuns}
\end{figure}

\fig{\ref{fig.EarthModelEigs}} shows a few approximate eigenvalues of the quadratic eigenproblem~\eqref{eq.HcurlQuadEig} arising in the electromagnetic wave propagation through our three-layer Earth model. For the sake of simplicity and to have a better visual representation of eigenvalues, we assume ${\mu=1}$\., and ${\eps=1}$\..
We represent the approximate eigenfunctions associated with two arbitrary eigenmodes in \fig{\ref{fig.EarthModelEigFuns}} to see how the electromagnetic wave looks like in different layers.
The electrical conductivity (i.e., the inverse of resistivity) performs as a {\em proportional} damping. Thus, eigenfunctions show a non-decaying oscillatory pattern within each layer. The decaying effect of material conductivities, however, appears in the time response that is not studied here (see~\cite{Cooke2000} for further details).

\vspace{-10pt}

\subsubsection{Acoustic fluid in a cavity with one absorbing wall}
\label{ssub:AcousticFluid}

\begin{figure}[!h]
	\centering
	\includegraphics[width=0.45\textwidth]{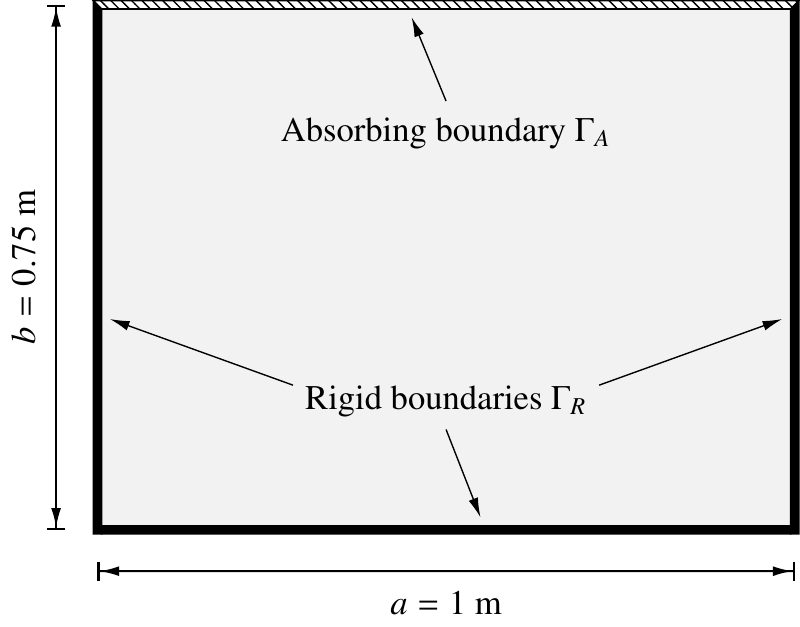}
	\caption{Physical domain of the vibroacoustic eigenproblem. We assume the cavity is filled with air and has one absorbing wall.}
	\label{fig.AcousticFluid}
\end{figure}

A practical application of the acoustic fluid in a cavity with absorbing walls occurs in the problem of decreasing noise level in, for example, a vehicle cabin. Absorbing walls, commonly covered by a thin layer of a viscoelastic material, can dissipate the fluid's acoustic energy.  In here, we consider a 2D rectangular cavity with one absorbing wall, as~\fig{\ref{fig.AcousticFluid}} depicts. We assume the cavity is filled with air ${(\,\rho=1~{\rm kg/m^3})}$\., the acoustic speed is ${c=340~{\rm m/s}}$\., and the viscoelastic material has the impedance parameters ${\alpha=5\times10^4~{\rm N/m^3}}$ and ${\beta=200~{\rm N\cdot s/m^3}}$.
\vspace{5pt}

\vspace{3pt}

\begin{figure}[!h]
	\centering
	\includegraphics{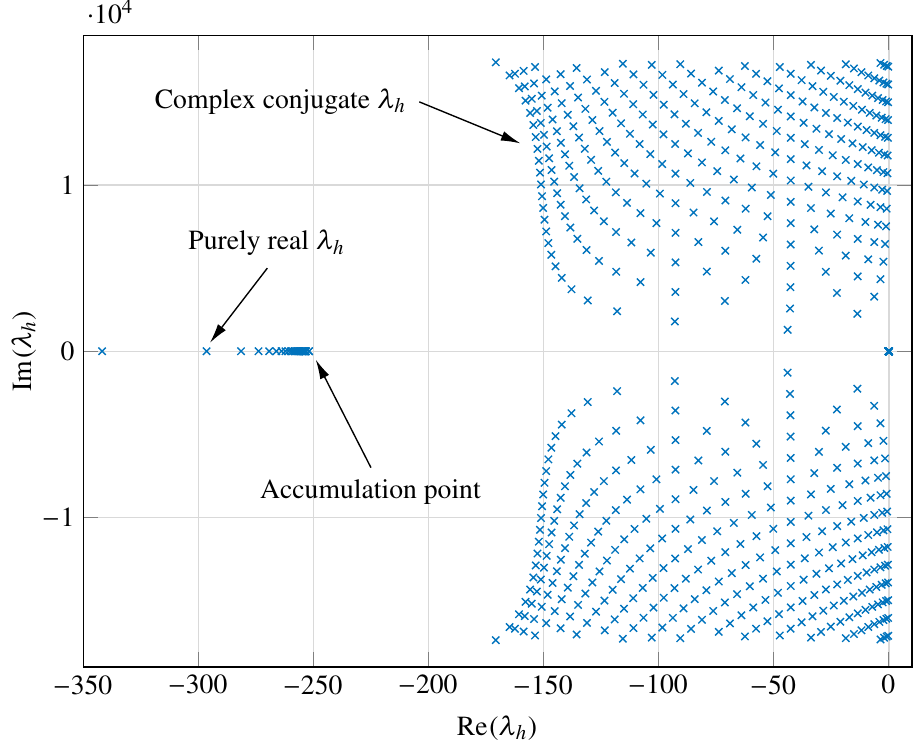}
	\caption{Approximate eigenvalues for quadratic eigenproblem arising from the acoustic fluid motion in a rectangular cavity with one absorbing wall; we represent only the first accumulation point of purely real eigenvalues at ${-\frac{\alpha}{\beta}=-250~{\rm rad/s}}$.}
	\label{fig.AcousticFluidEigs}
\end{figure}

\begin{figure}[!h]
	\begin{subfigure}{0.5\textwidth}\hspace{-0.2cm}
		\begin{subfigure}{0.24\textwidth}\hspace{-0.5cm}
			\includegraphics{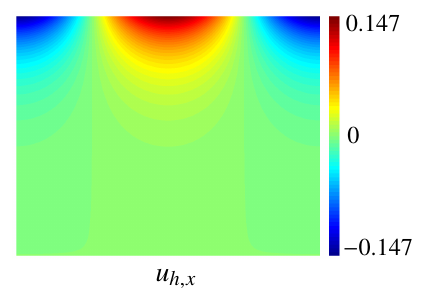}
		\end{subfigure}\hspace{2.07cm}
		\begin{subfigure}{0.24\textwidth}\centering
			\includegraphics{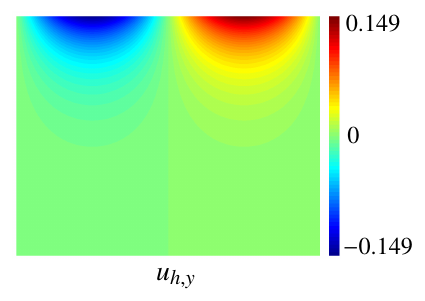}
		\end{subfigure}\hspace{3pt}
		\captionsetup{justification=centering}
		\caption{Numerical eigenfunctions associated with $\lmh=-341.80~{\rm rad/s}$}
	\end{subfigure}
	\begin{subfigure}{0.5\textwidth}\hspace{-0.15cm}
		\begin{subfigure}{0.24\textwidth}
			\includegraphics{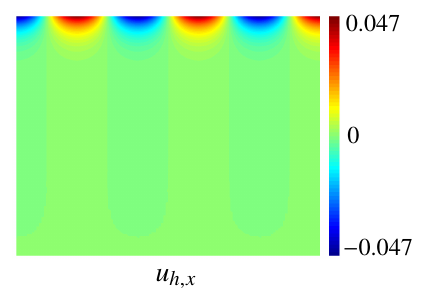}
		\end{subfigure}\hspace{2.07cm}		
		\begin{subfigure}{0.24\textwidth}\centering
			\includegraphics{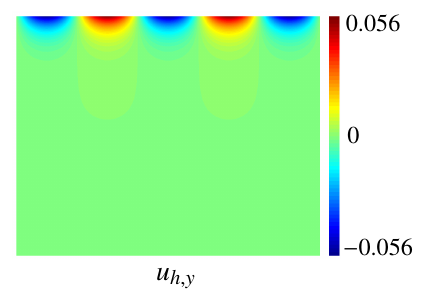}
		\end{subfigure}
		\captionsetup{justification=centering}
		\caption{Numerical eigenfunctions associated with $\lmh=-3480.48~{\rm rad/s}$}
	\end{subfigure}
	\caption{Two arbitrary normalized eigenfunctions corresponding to two distinct purely real eigenvalues in two different branches of the vibroacoustic case study.}
	\label{fig.AcousticFluidEigFuns}
\end{figure}

\fig{\ref{fig.AcousticFluidEigs}} shows a few approximate eigenvalues of the second case study. Again, the eigenvalues form a combination of conjugate complex pairs and purely real values.  The negative real part of the eigenvalues connotes the decaying behavior of vibration modes of the fluid. Purely real eigenvalues, corresponding to overdamped modes, theoretically have two {\em accumulation} points at $-\frac{\alpha}{\beta}$ and $-\frac{2\alpha}{\beta}$ (see~\cite{ Bermdez1999, Bermdez2000} for more details). More precisely, the eigenvalues exist in two separate branches, namely ${-\frac{2\alpha}{\beta}<\lmh<-\frac{\alpha}{\beta}}$ and ${\lmh<-\frac{2\alpha}{\beta}}$\,, which correspond to two possible solutions of the {\em dispersion} equations (see \Sec{\ref{sub:AccuracyTests}}).  \fig{\ref{fig.AcousticFluidEigFuns}} illustrates two eigenfunctions of two distinct eigenmodes in two different branches. The dissipation effect of the absorbing boundary causes the solution to decay in the vertical direction while it is oscillatory along the horizontal one.

\vspace{2pt}

\begin{remark}
  When the stiffness, mass, and damping matrices have widely different norms, which could be the case in the current vibroacoustic eigenproblem, the solution of the linearized eigensystem~\eqref{eq.Linearization} may suffer a significant numerical error (c.f.,~\cite{ Tisseur2000}). We overcome this issue, as recommended in~\cite{ Tisseur2000}, by solving a scaled eigenproblem:
  \begin{align}
  \big(\K+\gamma\.\C_\varsigma+\gamma^{\.2}\M_\varsigma\big)\.\tuu=\textbf{0}\,,
  \end{align}
  in place of the original one, where 
  ${\gamma=\lmh/\varsigma}$\,,
  ${\C_\varsigma=\varsigma\.\C}$\., and
  ${\M_\varsigma=\varsigma^{\.2}\M}$\.. In this setup, the scaling factor $\varsigma$ is
  \begin{align}
  \varsigma \ceq \sqrt{\dfrac{\norm{\K}_\infty}{\norm{\M}_\infty}}\,,
  \end{align}
  where the infinity norm is obtained as the maximum absolute value of respective matrices.
\end{remark}


\subsection{Computational efficiency tests}
\label{sub:EfficiencyTests}

We report the computational cost of the eigensolution of our two cases studies. Since the computational cost of solving quadratic eigensystems discretized by either \ibH(curl) or \ibH(div) spaces are approximately the same, we only present the costs of solving the eigenproblem of electromagnetic wave propagation.  We test different mesh sizes with ${\ne=2^s}$ ${(s=5,6,...,11)}$ elements in each direction, and different polynomial degrees of B-spline bases, namely ${\p=3,4,5,6}$.  We assess the computational efficiency of the rIGA framework by considering an optimal macroelement size for rIGA discretization; that is, each macroelement consists of a ${16\times16}$ grid (c.f.,~\cite{ Garcia2017, Garcia2019, Hashemian2021}).

\begin{figure}[!h]
	\begin{subfigure}{1\textwidth}\hspace{-0.2cm}
		\includegraphics{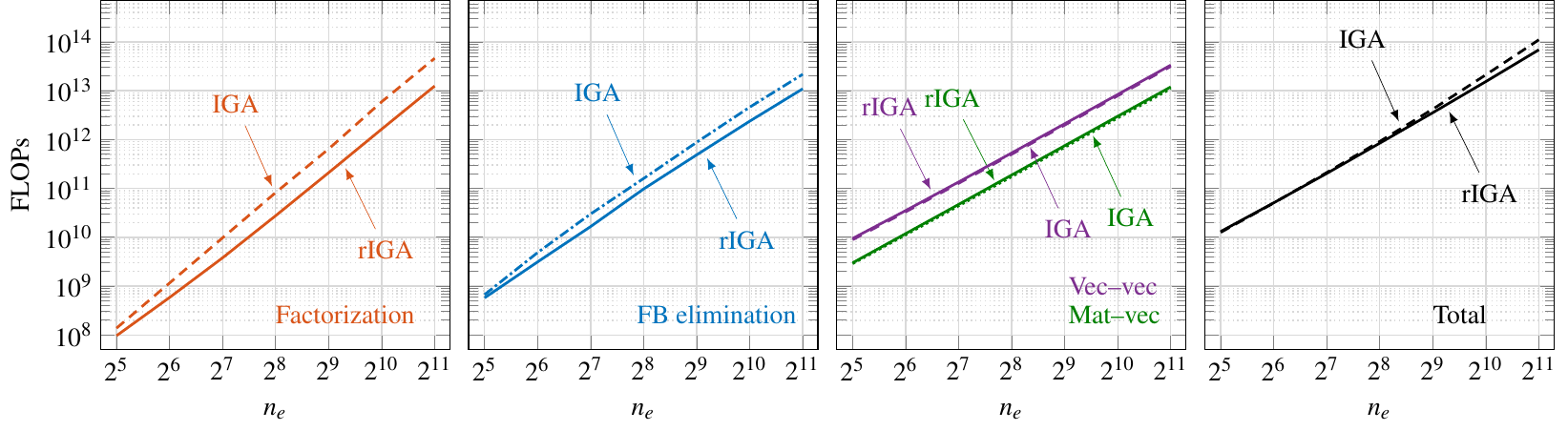}\vspace{-0.1cm}
		\caption{$\p=3$}\vspace{0.25cm}
	\end{subfigure}
	\begin{subfigure}{1\textwidth}\hspace{-0.2cm}
		\includegraphics{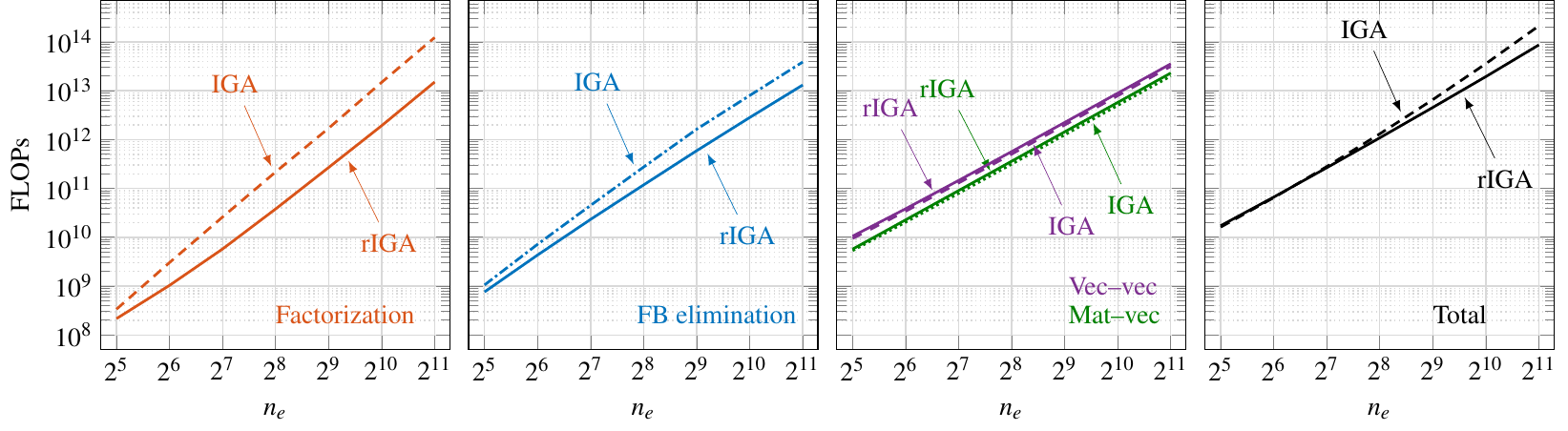}\vspace{-0.1cm}
		\caption{$\p=4$}\vspace{0.25cm}
	\end{subfigure}
	\begin{subfigure}{1\textwidth}\hspace{-0.2cm}
		\includegraphics{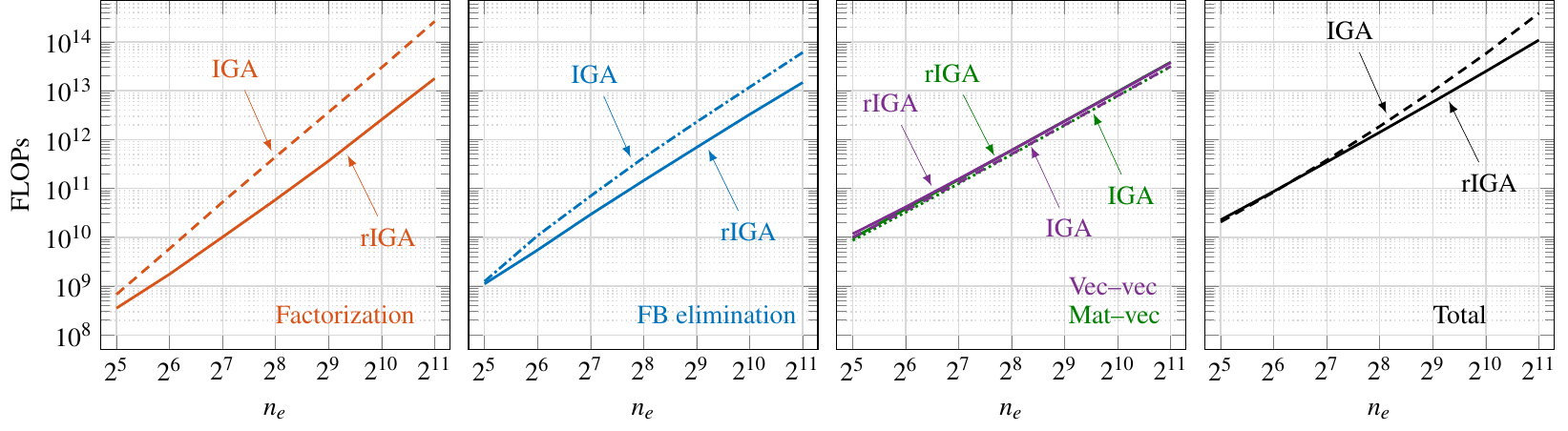}\vspace{-0.1cm}
		\caption{$\p=5$}\vspace{0.25cm}
	\end{subfigure}
	\begin{subfigure}{1\textwidth}\hspace{-0.2cm}
		\includegraphics{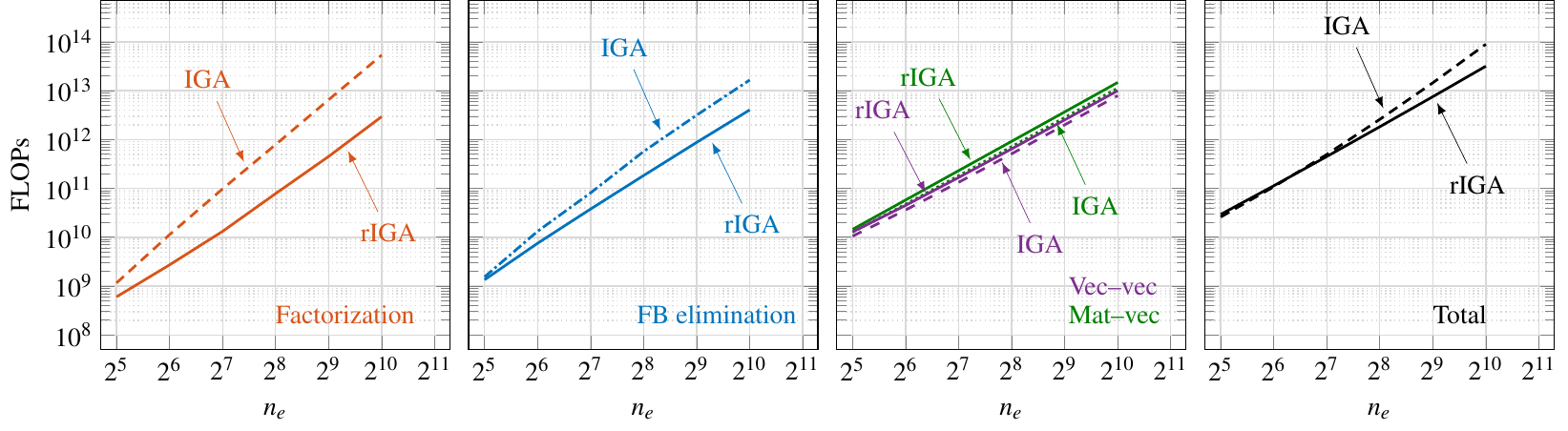}\vspace{-0.1cm}
		\caption{$\p=6$}
	\end{subfigure}
	\caption{Computational cost (FLOPs number) of different numerical operations employed by the eigensolver when using IGA and rIGA discretizations. We report the cost of finding 500 eigenpairs of quadratic eigenproblems discretized by \ibH(curl) and \ibH(div) spaces. Due to highly demanding computations, for ${\ne=2048}$, we only test up to ${\p=5}$.\\[-10pt]}
	\label{fig.FlopTest}
\end{figure}

\begin{figure}[!h]\centering
	\begin{subfigure}{1\textwidth}\centering
		\includegraphics{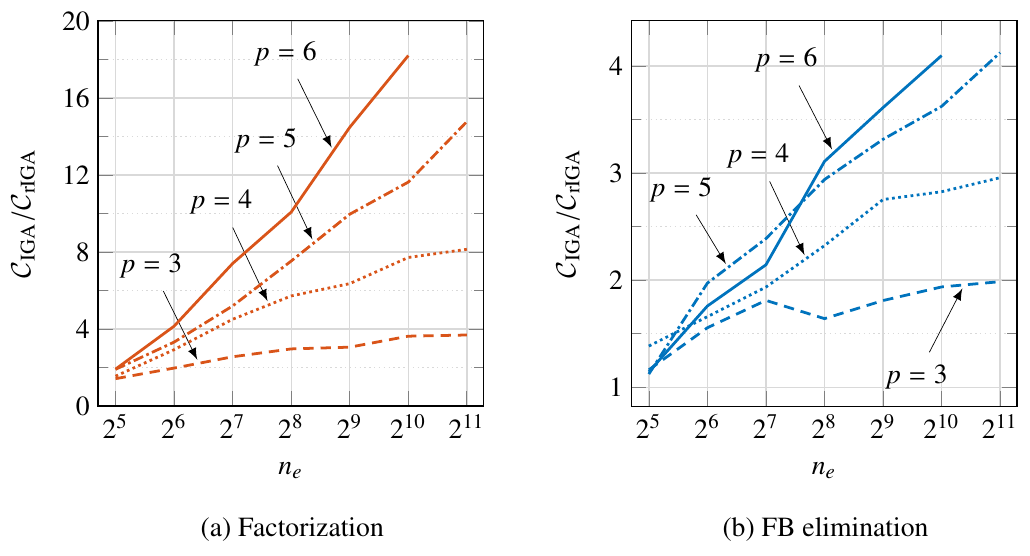}
	\end{subfigure}\vspace{5pt}
	\begin{subfigure}{1\textwidth}\centering
		\includegraphics{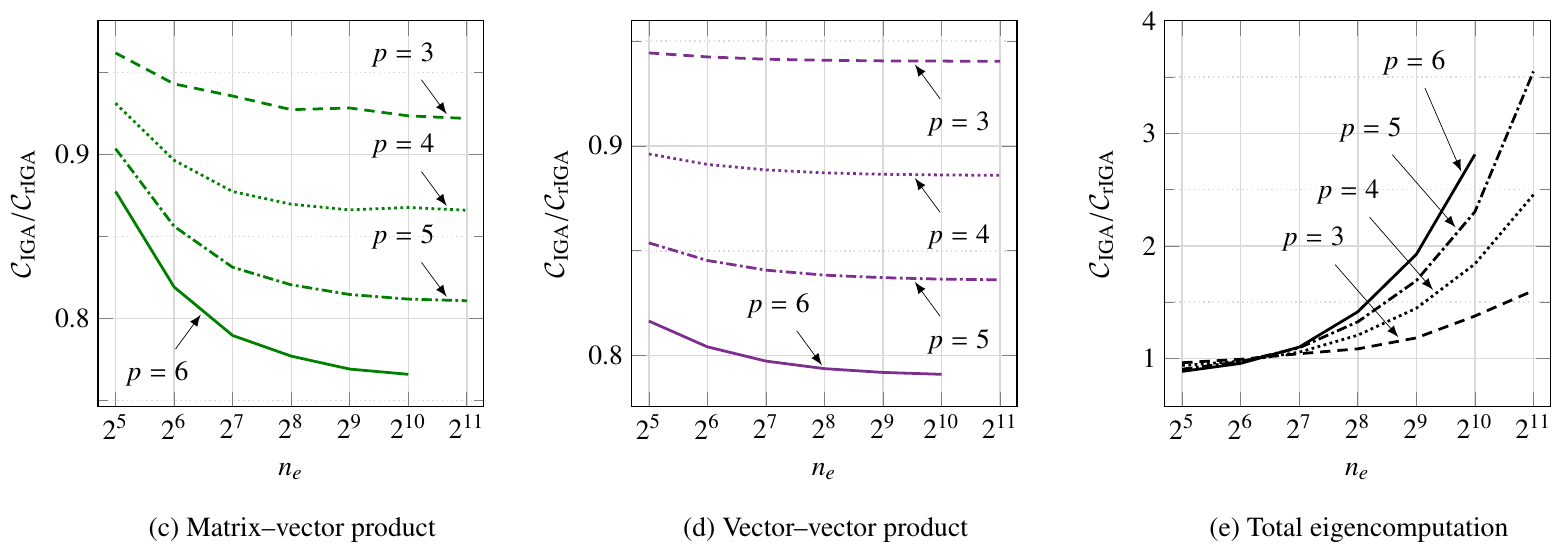}
	\end{subfigure}
	\caption{Cost improvement (or degradation) rates of different numerical operations of the eigenanalysis for rIGA when finding 500 eigenpairs of quadratic eigenproblems discretized by \ibH(curl) and \ibH(div) spaces. \vspace{12pt}}
	\label{fig.FlopTestImprovements}
\end{figure}

\fig{\ref{fig.FlopTest}} describes the contribution of each of the most expensive operations to the total eigencomputation cost when using maximum-continuity and refined IGA frameworks to approximate ${\Nev=500}$ eigenpairs.
We also investigate the cost improvement (or degradation) rates of each operation under the employment of rIGA by plotting the relative cost ${\Cr_\IGA/\Cr_\rIGA}$ against the mesh size in \fig{\ref{fig.FlopTestImprovements}}.
For large problems, results indicate an improvement in the factorization cost close to ${\Or\big((\p-1)^2\big)}$ and of approximately ${\Or(\p-1)}$ for forward/backward eliminations. There is a slight degradation in the cost of matrix--vector and vector--vector multiplications due to the slight increase of the system size in rIGA. 
In summary, the total observed cost saving for the entire eigensolution is up to ${\Or(\p-1)}$ for large domains (see \tab{\ref{tab:FlopTest}}). 
If the problem size is large, the total computational cost is governed only by matrix factorization. Therefore, we predict the total time improvements of up to ${\Or\big((\p-1)^2\big)}$. To observe this scaling, we would need larger computational resources. On the other hand, for small problems (e.g., systems with ${\ne\leq128}$), the overall cost of maximum-continuity IGA is comparable (or even smaller) to that of rIGA, which occurs when the matrix factorization cost is a small fraction of the total cost.  Finally, the number of requested eigenpairs affects our cost improvements (see~\Sec{\ref{sub:CostSummary}}). However, in the case of a very large $\Nev$\., the eigenvalues may be located far from the target shift, thus, decreasing the accuracy of the spectral approximation. Consequently, a reasonable number of eigenpairs should be sought.

\begin{table}[!h]
	\centering
	\caption{Improvement/degradation rates of different numerical operations of the eigenanalysis when using rIGA discretizations. We report the results associated with ${\ne=1024,2048}$ in \fig{\ref{fig.FlopTestImprovements}} to show the asymptotic improvement factor of rIGA.}
	\label{tab:FlopTest}
	\small
	\begin{tabular}{@{}lp{1.25cm}lllll@{}}
		\toprule
		        				&       & Factorization & FB elimination & Mat--vec    & Vec--vec    & Total cost \\
    	Degree  				& $\ne$ & improved by   & improved by    & degraded by & degraded by & improved by \\[1.5pt]\toprule
		\multirow{2}{*}{$\p=3$}	& 1024  & 3.634         & 1.937          & 0.923       & 0.940       & 1.376 \\
						        & 2048  & 3.696         & 1.986          & 0.921       & 0.940       & 1.602 \\ \midrule
		\multirow{2}{*}{$\p=4$} & 1024  & 7.721         & 2.825          & 0.867       & 0.886       & 1.839 \\
						        & 2048  & 8.145         & 2.957          & 0.865       & 0.886       & 2.454 \\ \midrule
		\multirow{2}{*}{$\p=5$} & 1024  & 11.634        & 3.622          & 0.811       & 0.836       & 2.303 \\
						        & 2048  & 14.775        & 4.124          & 0.810       & 0.836       & 3.549 \\ \midrule
		\multirow{2}{*}{$\p=6$} & 1024  & 18.203        & 4.097          & 0.765       & 0.791       & 2.811 \\ 
            					& 2048  & ***           & ***            & ***         & ***         & ***   \\
    \bottomrule
	\end{tabular}
\end{table}


\subsection{Accuracy versus cost}
\label{sub:AccuracyTests}

We present the results of accuracy-versus-cost tests when using different IGA and rIGA discretizations to solve eigenproblems arising in electromagnetics and vibroacoustics.  In~\cite{ Puzyrev2018, Hashemian2021}, the authors show for linear Hermitian definite pencils that, other than for the outliers, both IGA and rIGA frameworks result in almost the same accuracy in the eigensolution.  In the following, we consider one arbitrary eigenmode for each case study and compare the eigenvalue and eigenfunction errors. We show that rIGA delivers almost the same accuracy than IGA, but with a lower computational effort.


\subsubsection{Electromagnetic wave propagation eigenproblem}
\label{ssub:EarthModelAccuracy}

The analytic eigensolution of~\eqref{eq.HcurlStrong1} in conductive media is unavailable. Thus, we consider ${\pmb{\sg}=\textbf{0}}$ and study the accuracy of approximate eigensolution. The analytic eigenvalues and eigenfunctions of the electromagnetic wave in a non-conductive unit square are expressed in a tensor form as follows (see, e.g.,~\cite{BUFFA2010}):
\begin{align}
\lm_{ij}&=\pi^{\.2}(i^{\.2}+j^{\.2})\,,\\
\EE_{ij}&=\dfrac{2}{\sqrt{i^{\.2}+j^{\.2}}}\begin{bmatrix}
-j\cos(i\pi x)\sin(j\pi y)\\[2pt]
~i\sin(i\pi x)\cos(j\pi y)
\end{bmatrix}\..
\label{eq.ExactSolElectromagnetic}
\end{align}

\begin{figure}[!h]
	\begin{subfigure}{0.49\textwidth}\hspace{-0.25cm}
		\includegraphics{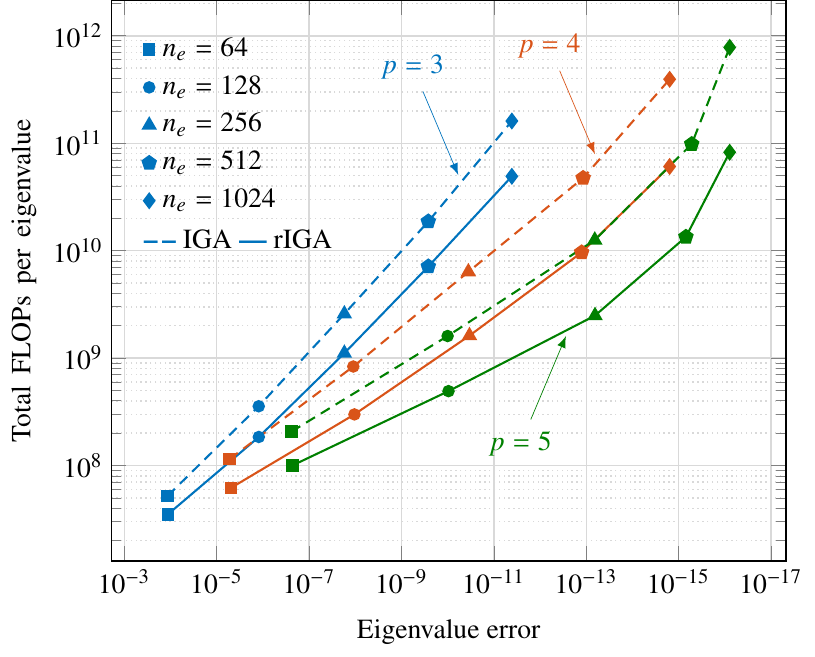}
		\caption{Eigenvalue error: $(\lmh-\lm)/\lm$}
		\label{fig.AccuracyTestElectromag.a}
	\end{subfigure}
	\begin{subfigure}{0.49\textwidth}\hspace{0.17cm}
		\includegraphics{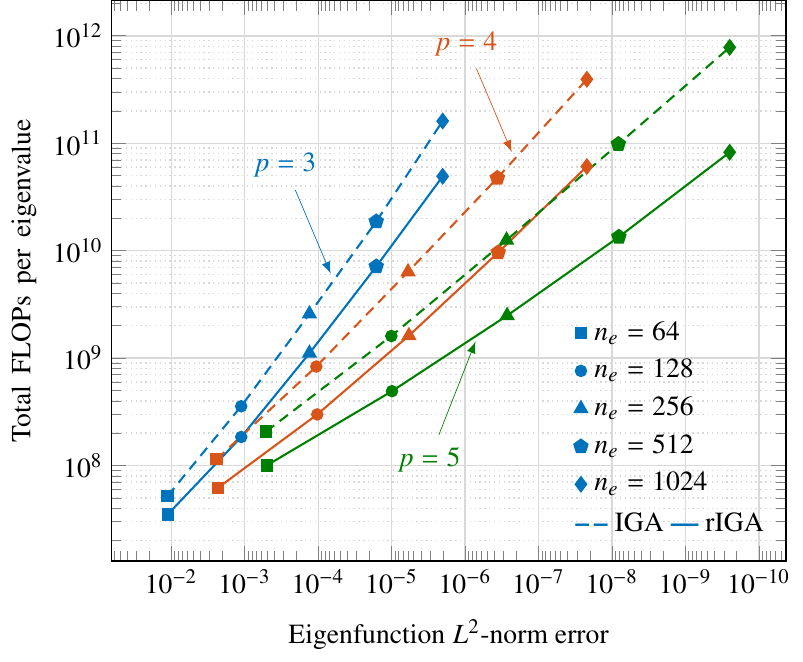}
		\caption{Eigenfunction $L^2$-norm error: $\normL{\.\EEh-\EE\.}$}
		\label{fig.AccuracyTestElectromag.b}
	\end{subfigure}
	\caption{Accuracy versus eigencomputation cost for the electromagnetic eigenvalue problem in a non-conductive unit square. We compare the results of different IGA and rIGA discretizations when approximating the eigenpair ${(\lm_{23,23},\EE_{23,23})}$.}
	\label{fig.AccuracyTestElectromag}
\end{figure}

\figs{\ref{fig.AccuracyTestElectromag.a}}{\ref{fig.AccuracyTestElectromag.b}} show the total computational cost per eigenmode versus the eigenvalue error ${(\lmh-\lm)/\lm}$ and eigenfunction $L^2$-norm error ${\normL{\.\EEh-\EE\.}}$\., respectively. We consider the analytic eigenvalue ${\lm_{23,23}=10442.04~{\rm rad/s}}$ and use IGA and rIGA discretizations with different polynomial degrees and mesh sizes.
Both figures confirm that under the employment of rIGA, we obtain almost the same approximation error as the maximum-continuity IGA with a lower computational effort.
Here, in the absence of conductivity, the quadratic eigenproblem is converted to a definite generalized one. Thus, we follow the notations described in \Rem{\ref{rem.LinearHcurl}} and use $\lm$ in reference to eigenvalues of the electromagnetic wave. Additionally, we use the generalized eigensolver to solve the eigenproblem, thus, obtaining higher improvements in the number of FLOPs when using rIGA.
\vspace{5pt}

\subsubsection{Acoustic fluid eigenproblem}
\label{ssub:AcousticFluidAccuracy}

\begin{figure}[!h]
	\begin{subfigure}{0.49\textwidth}\hspace{-0.25cm}
		\includegraphics{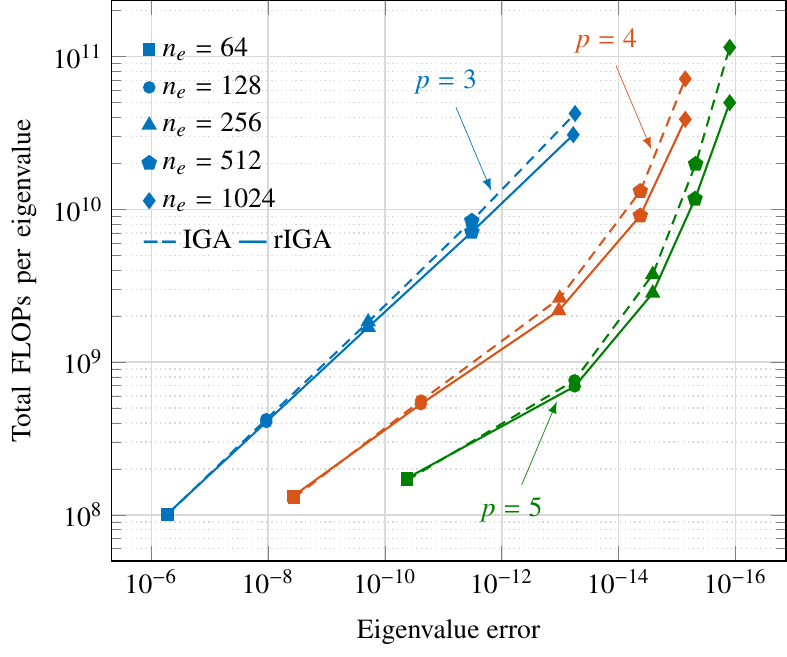}
		\caption{Eigenvalue error: $(\lmh-\lm)/\lm$}
		\label{fig.AccuracyTestAcoustic.a}
	\end{subfigure}
	\begin{subfigure}{0.49\textwidth}\hspace{0.17cm}
		\includegraphics{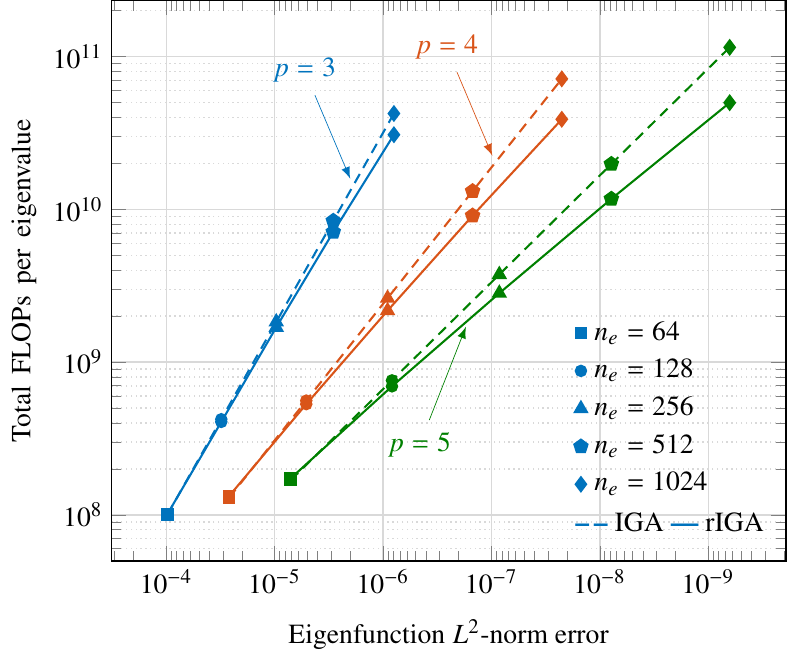}
		\caption{Eigenfunction $L^2$-norm error: $\normL{\.\uu_h-\uu\.}$}
		\label{fig.AccuracyTestAcoustic.b}
	\end{subfigure}
	\caption{Accuracy versus eigencomputation cost for the quadratic eigenproblem arising in vibroacoustics. We compare the results of different IGA and rIGA discretizations when approximating the eigenpair ${(\lm_{10},\uu_{10})}$ on the second eigenbranch of the rectangular cavity of \fig{\ref{fig.AcousticFluid}}.}
	\label{fig.AccuracyTestAcoustic}
\end{figure}

We obtain the analytical vibration modes of the acoustic fluid inside a 2D rectangular cavity with one absorbing wall (\fig{\ref{fig.AcousticFluid}}) by separation of variables as follows (see~\cite{Bermdez1999}):
\begin{align}
\uu_j=\dfrac{-1}{\rho\lm_j^2}\nabla\pra{\cos\Big(\dfrac{j\pi x}{a}\Big)\cosh(\eta_jy)}=
\dfrac{-1}{\rho\lm_j^2}\begin{bmatrix}
-\dfrac{j\pi}{a}\sin\Big(\dfrac{j\pi x}{a}\Big)\cosh(\eta_jy)\\[8pt]
~\eta_j\cos\Big(\dfrac{j\pi x}{a}\Big)\sinh(\eta_jy)
\end{bmatrix}\.,
\label{eq.ExactSolAcoustic}
\end{align}
where eigenvalues $\lm_j$ and auxiliary variables ${\eta_j\in\CC}$ are solutions of the dispersion equations:
\begin{align}
\eta_j^{\.2} &= \dfrac{\lm_j^2}{c^{\.2}}+\dfrac{\pi^{\.2}}{a^{\.2}}j^{\.2}\,, 
\label{eq.dispersion.a}\\[3pt]
\eta_j\tanh(\eta_j\.b) &= -\dfrac{\rho\lm_j^2}{\alpha+\lm_j\.\beta}\,.
\label{eq.dispersion.b}
\end{align}
Herein, we seek purely real eigenvalues.  As~\Sec{\ref{ssub:AcousticFluid}} states, for each $j$, the eigenvalues have two branches that accumulate at $-\frac{\alpha}{\beta}$ and $-\frac{2\alpha}{\beta}$.  \figs{\ref{fig.AccuracyTestAcoustic.a}}{\ref{fig.AccuracyTestAcoustic.b}} show the eigenvalue error ${(\lmh-\lm)/\lm}$ and eigenfunction $L^2$-norm error ${\normL{\.\uuh-\uu\.}}$ versus the total eigencomputation cost per eigenmode when using different IGA and rIGA discretizations. The results are associated with ${\lm_{10}=-7377.39~{\rm rad/s}}$ on the second branch (i.e., ${\lm<-\frac{2\alpha}{\beta}}$).  Again, the plots confirm that the implementation of rIGA reduces the computational cost for almost the same numerical error.

\section{Conclusions} 
\label{sec.Conclusions}

We use refined isogeometric analysis (rIGA) discretizations to solve quadratic eigenvalue problems (QEPs). We compare the computational performance of rIGA versus that of maximum-continuity IGA when solving quadratic eigensystems arising in 2D multifield problems discretized by \ibH(curl) and \ibH(div) spaces.  For large problem sizes, the most expensive numerical operation of the eigenanalysis is matrix factorization, followed by matrix--vector and vector--vector operations in the sense of Krylov projections. We estimate that the LU factorization is asymptotically ${\Or\big((\p-1)^2\big)}$ times faster with an rIGA discretization compared to IGA.  As a result, we can theoretically reach an improvement of ${\Or\big((\p-1)^2\big)}$ in the total eigencomputation cost of vector-valued multifield eigenproblems ---\.while it is expected to be ${\Or(\p^2)}$ in scalar-valued problems in $H^1$.  Our numerical tests show that the computational savings associated with the implementation of rIGA tend to ${\Or(\p-1)}$, as we need to analyze larger systems to observe the asymptotic behavior.  For small problems, the improvement hardly reaches the expected rates. Thus, we suggest using the maximum-continuity IGA discretization only for small problems.  Finally, quadratic eigensystems have complex eigenvalues in most occasions, forcing the eigensolver to project the problem onto a larger Krylov subspace (as compared to the definite generalized eigenproblems). Thus, more matrix--vector and vector--vector operations are required, resulting in a deterioration of the improvement rate of rIGA versus IGA when solving quadratic eigenproblems.


\section*{Acknowledgment}
This work has received funding from
the European Union's Horizon 2020 research and innovation program under the Marie~Sklodowska-Curie grant agreement No.~777778 (MATHROCKS); the European Regional Development Fund (ERDF) through the Interreg V-A Spain-France-Andorra program POCTEFA 2014-2020 Project PIXIL (EFA362/19); the Spanish Ministry of Science and Innovation projects with references PID2019-108111RB-I00 (FEDER/AEI) and PDC2021-121093-I00, the “BCAM Severo Ochoa” accreditation of excellence (SEV-2017-0718); and the Basque Government through the BERC 2018-2021 program, the three Elkartek projects 3KIA (KK-2020/00049), EXPERTIA (KK-2021/00048), and SIGZE (KK-2021/00095), and the Consolidated Research Group MATHMODE (IT1294-19) given by the Department of Education.
This publication was also made possible in part by the Professorial Chair in Computational Geoscience at Curtin University and the Deep Earth Imaging Enterprise Future Science Platforms of the Commonwealth Scientific Industrial Research Organisation, CSIRO, of Australia. The Curtin Corrosion Centre and the Curtin Institute for Computation kindly provide ongoing support.
The authors also acknowledge the computer resources and technical support provided by Barcelona Supercomputing Center through the MareNostrum4 cluster (activity IDs IM-2020-3-0009 and IM-2021-2-0015).


\section*{References}
\bibliography{_eig_refs}

\end{document}